\documentclass[leqno]{aoscmsize11}%
\usepackage{amsfonts}
\usepackage{amssymb}
\usepackage{amsthm}
\usepackage{latexsym}
\usepackage{graphicx}
\usepackage{amsmath}%
\newtheorem{theorem}{Theorem}

\newtheorem{lemma}{Lemma}

\def\theequation{\arabic{equation}}
\NONUMBIB
\def\baselinestretch{1.5}
\begin{document}

\SPECFNSYMBOL{}{1}{}{}{}{}{}{}{}
%
\AOSMAKETITLE
%
%

\AOSAMS{Primary 62F12, 62M05; secondary 60H10, 60J60.}
\AOSKeywords{{Jumps, efficiency, inference, discrete sampling}\\}
\AOStitle{FISHER'S INFORMATION FOR DISCRETELY SAMPLED LEVY PROCESSES}
\AOSauthor{Yacine A\"{\i}t-Sahalia\footnote
{Supported in part by NSF Grants SBR-0111140 and SBR-0350772.}
and Jean Jacod}
\AOSaffil{Princeton University and Universit\'e de Paris-6}
\AOSlrh{YACINE AIT-SAHALIA AND JEAN JACOD}
\AOSrrh{INFORMATION FOR LEVY PROCESSES}
\renewcommand{\baselinestretch}{1.0}
\AOSAbstract
{This paper studies the asymptotic behavior of the Fisher information for a
L{\'{e}}vy process discretely sampled at an increasing frequency.
We show that it is possible to distinguish not only the continuous part of the process from its jumps part,
but also different types of jumps, and derive the rates of convergence of efficient estimators.}
\maketitle
%

\BACKTONORMALFOOTNOTE{1}
%
%

\renewcommand{\baselinestretch}{1.5}%

\section{Introduction.\label{sec:intro}}

Models allowing for sample path discontinuities are of considerable interest
in mathematical finance, for instance in option pricing [see e.g.,
Eberlein and Jacod (1997), Chan (1999), Boyarchenko and Levendorskii (2002),
Mordecki (2002) and Carr and Wu (2004)], testing for the presence of jumps in
asset prices [see A\"\i t--Sahalia (2002) and Carr and Wu (2003)],
interest rate modelling [see e.g., Eeberlein abd Raible (1999)], risk
management [see e.g., Eberlein et al. (1998) and Khindanova et
al. (2001)], optimal portfolio choice [see e.g., Kallsen (2000),
Rachev and Han (2000) and Emmer and Kl\"uppelberg (2004)], stochastic
volatility modelling [see e.g., Barndorff--Nielsen (1997),
Barndorff--Nielsen (1998), Leblanc and Yor (1998), Carr et al. (2003)
and Kl\"uppelberg et al. (2004)] or for the
purpose of better describing asset returns data [see e.g., Mandelbrot
(1963), Fama and Roll (1965), Mittnik and Rachev (2001), Carr et al. (2002)].

While these theoretical models are commonly used in mathematical finance,
relatively little is known about the corresponding inference problem, which is
a difficult one. A string of the literature focuses on the tail properties of
stable processes to estimate the stable index [see e.g., Fama and Roll
(1968), Fama and Roll (1971), de Haan and Resnick (1980), Dumouchel (1983),
McCulloch (1997)]. Since L\'{e}vy processes have known characteristic
functions, given by the L\'{e}vy-Khintchine formula, a method often proposed
is based on the empirical characteristic function as an estimating equation
[see e.g., Press (1972), Fenech(1976), Feuerverger and McDunnough (1981b),
Chapter 4 in Zolotarev (1986) and Singleton (2001)], maximum likelihood by
Fourier inversion of the characteristic function [see
Feuerverger and McDunnough (1981a)], or a regression based on the explicit form
of the characteristic function [see Koutrouvelis (1980)]. Some of these
methods were compared in Akgiray and Lamoureux (1989).

Fairly little is known in most cases as to the optimality of statistical
procedures in the presence of jumps. So we consider in this paper the behavior
of the Fisher information when the observations are generated by a L\'{e}vy
process $X$ whose law depends on a parameter vector $\eta$ to be estimated. In
light of the Cramer-Rao bound, our objective is to establish the optimality of
potential estimators of $\eta$, and the rate at which they will converge.
While we focus on its implications for the classical likelihood inference
problem, Fisher's information also plays as usual an important role in
Bayesian inference or in determining the form of asymptotically most
powerful tests.

The essential difficulty in this class of problems is the fact that the
density of most discretely sampled L\'{e}vy processes, hence the corresponding
likelihood function and Fisher's information, are not known in closed form.
Representations in terms of special functions are available [see
Zolotarev~(1995) in terms of Meijer $G-$functions and
Hoffmann--J\o rgensen (1993) in terms of incomplete hypergeometric functions]
although they do not appear to lead to practical formulae. One must therefore
rely on numerical methods as the maximum likelihood estimator cannot be
computed exactly [see Dumouchel (1971) for a multinomial approximation to
the likelihood function, and Nolan~(1997) and Nolan (2001)]. Therefore,
there is potential value in considering alternative estimators which can both
be computed explicitly and be rate-efficient. Indeed, in a companion paper
[A\"\i t--Sahalia and Jacod (2004)], we propose estimators designed to
achieve the efficient
rate that we identify in this paper based on the convergence properties of the
Fisher information.

Let us be more specific. The L\'{e}vy process $X$ is observed at $n$ times
$\Delta,2\Delta,\ldots n\Delta$. Recalling that $X_{0}=0$, this amounts to
observing the $n$ increments $X_{i\Delta}-X_{(i-1)\Delta}$. So when $\Delta>0$
is fixed, we observe $n$ i.i.d.\thinspace\ variables distributed as
$X_{\Delta}-X_{0}$ and having a density which depends smoothly on the
parameter $\eta$, and we are on known grounds: the Fisher information at stage
$n$ has the form $I_{n,\Delta}(\eta)=nI_{\Delta}(\eta)$, where $I_{\Delta
}(\eta)>0$ is the Fisher information (an invertible matrix if $\eta$ is
multi--dimensional) of the model based upon the observation of the single
variable $X_{\Delta}-X_{0}$; we have the LAN property with rate $\sqrt{n}$;
the asymptotically efficient estimators $\widehat{\eta}_{n}$ are those for
which $\sqrt{n}(\widehat{\eta}_{n}-\eta)$ converges in law to the normal
distribution $\mathcal{N}(0,I_{\Delta}(\eta)^{-1})$, and the MLE does the job
[see e.g., Dumouchel (1973a)].

Things become more complicated when the time interval separating successive
observations, $\Delta$, varies, and more specifically in the limit when it
becomes small. This type of asymptotics corresponds to a situation which is
increasingly common in financial applications, where high frequency data are
available with sampling intervals measured in seconds for many stocks or
currencies. At stage $n$ we have $n$ observations, recorded at times
$i\Delta_{n}$ for some time lag $\Delta_{n}$ going to $0$. Equivalently, we
observe the $n$ increments $\chi_{i}^{n}=X_{i\Delta_{n}}-X_{(i-1)\Delta_{n}}$
which are i.i.d.\ for any given $n$. These increments are known as the
log-returns when $X$ is the log of the price of a financial asset. However the
law of $\chi_{1}^{n}$ depends on $n$, and indeed weakly converges to the Dirac
mass at $0$. The Fisher information at stage $n$ still has the form
$I_{n,\Delta_{n}}(\eta)=nI_{\Delta_{n}}(\eta)$, but the behavior of the
information $I_{\Delta_{n}}(\eta)$ is far from obvious.

In order to say more about the behavior of Fisher's information, we need of
course to parametrize the model. For the same reason that the computation of
the MLE\ is hindered by the absence of an explicit density, the analysis of
the Fisher information matrix is difficult. Dumouchel (1973b) and
Dumouchel (1975) computed the information by numerical approximation of the
densities and their derivatives. Such direct computation is numerically
cumbersome because the series expansion for the density converges slowly,
especially when the order of the stable process is near one.
Brockwell and Brown (1980) propose an alternative numerical computation of the
information based on a Fourier series for the derivatives of the density.

In this paper, we are able to explicitly describe, in closed form, the
limiting behavior of the Fisher information when we restrict attention to a
special kind of L\'{e}vy process that is relevant to applications in financial
statistics. While our form of the process is undoubtedly restrictive, it is
nevertheless sufficiently rich to exhibit a surprising range of different
asymptotic behaviors for the Fisher information. In fact, we will show that
different rates of convergence are achieved for different parameters, and for
different types of L\'{e}vy processes. Rates depart from the standard
$\sqrt{n}$ in a number of different and often unexpected ways.

Specifically, we split $X$ into the sum of two independent L\'{e}vy processes,
with possibly one or two scale parameters. That is, we suppose that
\begin{equation}
X_{t}=\sigma W_{t}+\theta Y_{t}. \label{eq:X=sW+aY}%
\end{equation}
Here, we have $\sigma>0$ and $\theta\in\mathbb{R}$, and $W$ is a standard
symmetric stable process with index $\beta\in(0,2]$, and we are often
interested in the situation where $\beta=2$ and so $W$ is a Wiener process
(hence the notation used). As for $Y$, it is another L\'{e}vy process, viewed
as a perturbation of $W.$ In some applications, $Y$ may represent frictions
that are due to the mechanics of the trading process, or in the case of
compound Poisson jumps it may represent the infrequent arrival of relevant
information related to the asset. In the latter case, $W$ is then the driving
process for the ordinary fluctuations of the asset value. $Y$ is independent
of $W$, and its law is either known or is a nuisance parameter, and is
\emph{dominated}\ by $W$ in a sense stated below. For example, when $W$ is a
Brownian motion, this just means that $Y$ has no Brownian part; when $\beta
<2$, then $Y$ could for example be another stable process with index
$\alpha<\beta$, or a compound Poisson process. The parameter vector we then
consider is $\eta=(\sigma,\beta,\theta).$

If $Y$ is viewed as a perturbation of $W,$ then our interest in studying the
Fisher information lies in deciding whether we can estimate the parameter
$\sigma$, and also in some cases the index $\beta$ (the only two parameters on
which the law of the process $\sigma W$ depends) with the same degree of
accuracy as when the process $Y$ is absent, at least asymptotically. The
answer to this question is \textquotedblleft yes\textquotedblright. When $W$
is a Brownian motion this means that one can distinguish between the jumps due
to $Y$ and the continuous part of $X$, and this fact was already known in the
specific example of a Brownian motion coupled with either a Poisson or Cauchy
process [see A\"\i t--Sahalia (2004)]. It comes more as a surprise
when $\beta<2$: we
can then discriminate between the jumps due to $W$ and those due to $Y$,
despite the fact that both processes jump and we only have discrete observations.

The paper is organized as follows. In Section \ref{sec:setup}, we set up the
problem and define in particular the class of processes $Y$ that are dominated
by $W$. In Section \ref{sec:base}, we study the baseline case where
$X_{t}=\sigma W_{t}$ and establish the properties of the Fisher information in
the absence of the perturbation process $Y.$ In Section
\ref{sec:semiparametric}, we characterize the set of processes $Y$ whose
presence does not affect the estimation of the base parameters $(\sigma
,\beta)$. Then we study in Section \ref{sec:examples} the estimation problem
for the dominated scale parameter $\theta$. In this case, the results vary
substantially according to the structure of the process $Y$, and we illustrate
the versatility of the situation by displaying the variety of convergence
rates that arise.

Finally, we also briefly consider in Section \ref{sec:multiplicative} a
slightly different model, where
\begin{equation}
X_{t}=\sigma(W_{t}+Y_{t}). \label{eq:X=sW+sY}%
\end{equation}
Here it is natural to consider the law of $Y$ and the index $\beta$ as known,
and $\sigma$ to be the only parameter to be estimated. The results are again a
bit unexpected, namely one can do as well as when $Y$ is absent, and in some
instances (when the law of $Y$ is sufficiently singular) the presence of $Y$
can in fact help us improve the estimation of $\sigma$.

All proofs are in Section \ref{sec:proofs}.

\section{Setup.\label{sec:setup}}

The characteristic function of $W_{t}$ is
\begin{equation}
\mathbb{E}(e^{iuW_{t}})=e^{-t|u|^{\beta}/2} \label{eq:Wchar}%
\end{equation}
The factor $2$ above is unusual for stable processes when $\beta<2$, but we
put it here to ensure continuity between the stable and the Gaussian cases. As
is well known, when $\beta<2$ we have $\mathbb{E}(\left\vert W_{t}\right\vert
^{\rho})<\infty$ if and only if $0<\rho<\beta$, and the tails of $W_{1}$
behave according to $\mathbb{P}\left(  W_{1}>w\right)  \sim c_{\beta}/\beta
w^{\beta}$ as $w\rightarrow\infty$ (and symmetrically as $w\rightarrow-\infty
$), where the constants $c_{\beta}$ are given by
\begin{equation}
c_{\beta}=\left\{
\begin{array}
[c]{ll}%
\frac{\beta(1-\beta)}{4~\Gamma\left(  2-\beta\right)  ~\cos\left(  \beta
\pi/2\right)  } & \mbox{if }~\beta\neq1\\[2.4mm]%
\frac{1}{2\pi} & \mbox{if }~\beta=1
\end{array}
\right.  \label{eq:cbeta}%
\end{equation}
This follows from the series expansion of the density due to
Bergstr\o m (1952), the duality property of the stable densities of order
$\beta$ and $1/\beta$ [see e.g., Chapter 2 in Zolotarev (1986)], with an
adjustment factor to reflect our definition of the characteristic function in
(\ref{eq:Wchar}).

The law of $Y$ (as a process) is entirely specified by the law $G_{\Delta}$ of
the variable $Y_{\Delta}$ for any given $\Delta>0$. We write $G=G_{1}$, and we
recall that the characteristic function of $G_{\Delta}$ is given by the
L\'{e}vy-Khintchine formula
\begin{equation}
\mathbb{E}(e^{ivY_{\Delta}})=\exp\Delta\left(  ivb-{\frac{cv^{2}}{2}}+\int
F(dx)\left(  e^{ivx}-1-ivx1_{\{|x|\leq1\}}\right)  \right)  \label{eq:Ychar}%
\end{equation}
where $(b,c,F)$ is the \textquotedblleft characteristic
triple\textquotedblright\ of $G$ (or, of $Y$): $b\in\mathbb{R}$ is the drift
of $Y,$ and $c\geq0$ the local variance of the continuous part of $Y,$ and $F$
is the L\'{e}vy jump measure of $Y$, which satisfies $\int\left(  1\wedge
x^{2}\right)  F(dx)<\infty$ [see e.g., Chapter II.2 in
Jacod and Shiryaev (2003)].

The fact that $Y$ is \textquotedblleft dominated\textquotedblright\ by $W$ is
expressed by the property that $G$ belongs to the class $\mathcal{G}_{\beta}$
which we define as follows. Let first $\Phi$ be the class of all increasing
and bounded functions $\phi:~(0,1]\rightarrow\mathbf{R}_{+}$ having
$\lim_{x\downarrow0}\phi(x)=0$. Then we set
\begin{align}
\mathcal{G}(\phi,\alpha)~  &  =\mbox{the set of all infinitely divisible
distributions with ~$c=0$~and, for all} ~x\in(0,1],\label{eq:Gclass}\\
&  \left\{
\begin{array}
[c]{ll}%
x^{\alpha}F([-x,x]^{c})\leq\phi(x) \qquad & \mbox{if }~\alpha<2\\
x^{2}F([-x,x]^{c})\leq\phi(x)\quad\mbox{and}~~ \int_{\{|y|\leq x\}}%
|y|^{2}F(dy)\leq\phi(x)~\qquad & \mbox{if }~\alpha=2,
\end{array}
\right. \nonumber
\end{align}%
\begin{equation}
\mathcal{G}_{\alpha}=\cup_{\phi\in\Phi}~\mathcal{G}(\phi,\alpha).
\label{eq:BGclass}%
\end{equation}

We have $\lim_{x\downarrow0}x^{\alpha}F([-x,x]^{c})=0$ if and only if the
function $\phi(y)=\sup_{x\in(0,y]}x^{\alpha}F([-x,x]^{c})$ belongs to $\Phi$,
whereas $\int_{\{|y|\leq x\}}|y|^{2}F(dy)$ always decreases to $0$ as
$x\downarrow0$, so we also have another, simpler, description of
$\mathcal{G}_{\alpha}$ for all $\alpha\in(0,2]$:
\begin{equation}
\mathcal{G}_{\alpha}=\left\{  G~\text{is infinitely divisible},~c=0,~\lim
_{x\downarrow0}x^{\alpha}F([-x,x]^{c})=0\right\}  . \label{eq:BGclass2}%
\end{equation}
We also have for any $0<x<y\leq1$:
\[
x^{2}F([-x,x]^{c})\leq x^{2}F([-y,y]^{c})+\int_{\{\|z\|\leq y\}}z^{2}F(dz),
\]
from which we deduce that $x^{2}F([-x,x]^{c})\to$ as $x\downarrow0$ for any
infinitely divisible $G$. Therefore, $\mathcal{G}_{2}$ is indeed the set of
all infinitely divisible laws $G$ such that $c=0$. Obviously $\alpha
<\alpha^{\prime}$ implies $\mathcal{G}_{\alpha}\subset\mathcal{G}%
_{\alpha^{\prime}}$. If $G$ is a (non necessarily symmetric) stable law with
index $\gamma$ it belongs to $\mathcal{G}_{\alpha}$ for all $\alpha>\gamma$,
but not to $\mathcal{G}_{\gamma}$. If $Y$ is a compound Poisson process plus a
drift, then $G$ is in $\cup_{\alpha>0}\mathcal{G}_{\alpha}$. \smallskip

The variables under consideration have densities which depend smoothly on the
parameters, so Fisher's information is an appropriate tool for studying the
optimality of estimators. In the basic case of the model (\ref{eq:X=sW+aY}),
the law of the observed process $X$ depends on the three parameters $\sigma$,
$\beta$, $\theta$ to be estimated, plus on the law of $Y$ which is summarized
by $G$. The law of the variable $X_{\Delta}$ has a density which depends
smoothly on $\sigma$ and $\theta$, so that the $2\times2$ Fisher information
matrix (relative to $\sigma$ and $\theta$) of our experiment exists; it also
depends smoothly on $\beta$ when $\beta<2$, so in this case the $3\times3$
Fisher information matrix exists. In all cases we denote it by $I_{n,\Delta
_{n}}(\sigma,\beta,\theta,G)$, and it has the form
\[
I_{n,\Delta_{n}}(\sigma,\beta,\theta,G)=n~I_{\Delta_{n}}(\sigma,\beta
,\theta,G),
\]
where $I_{\Delta}(\sigma,\beta,\theta,G)$ is the Fisher information matrix
associated with the observation of a single variable $X_{\Delta}$. We denote
the elements of the matrix $I_{\Delta}(\sigma,\beta,\theta,G)$ as $I_{\Delta
}^{\sigma\sigma}(\sigma,\beta,\theta,G),$ $I_{\Delta}^{\sigma\beta}%
(\sigma,\beta,\theta,G),$ etc. ~ We may occasionally drop $G$, but at this
stage it is mentioned because it may appear as a nuisance parameter in our
model and we wish to have estimates for the Fisher information that are
\emph{uniform}\ in $G$, at least on some reasonable class of $G$'s. Let us
also mention that in many cases the parameter $\beta$ is indeed known: this is
particularly true when $W$ is a Brownian motion.

For the other model (\ref{eq:X=sW+sY}) the Fisher information for estimating
$\sigma$ (a positive number here) is now denoted by $I_{n,\Delta_{n}}^{\prime
}(\sigma,\beta,G)$, and it still has the form
\[
I_{n,\Delta_{n}}^{\prime}(\sigma,\beta,G)=n~I_{\Delta_{n}}^{\prime}%
(\sigma,\beta,G),
\]
with $I_{\Delta}^{\prime}(\sigma,\beta,G)$ being the Fisher information
associated with the observation of a single variable $X_{\Delta}$.

\section{The baseline case:\ estimating the parameters of the stable process
$X=\sigma W$.\label{sec:base}}

In this section we consider the base case $Y=0$, that is we observe the stable
process $X=\sigma W$ with scale parameter $\sigma>0$ and index parameter
$\beta\in(0,2]$. In our general framework, this corresponds to the situation
where $G=\delta_{0},$ a Dirac mass at $0,$ and we set the (now unidentified)
parameter $\theta$ to $0$, or for that matter any arbitrary value.

We have only the two parameters $\sigma$ and $\beta$ here, and our objective
in this section is to compute the Fisher information matrix in this case:
\[
I_{\Delta}(\sigma,\beta,0,\delta_{0})=\left(
\begin{array}
[c]{ccc}%
I_{\Delta}^{\sigma\sigma}(\sigma,\beta,0,\delta_{0}) & ~~ & I_{\Delta}%
^{\sigma\beta}(\sigma,\beta,0,\delta_{0})\\
I_{\Delta}^{\sigma\beta}(\sigma,\beta,0,\delta_{0}) &  & I_{\Delta}%
^{\beta\beta}(\sigma,\beta,0,\delta_{0})
\end{array}
\right)  .
\]
In future sections, we will examine how the terms in $I_{\Delta}(\sigma
,\beta,\theta,G)$ relate to those in $I_{\Delta}(\sigma,\beta,0,\delta_{0}).$

\subsection{The scale parameter $\sigma.$}

By the scaling property of symmetric stable processes, which says that
$W_{\Delta}$ and $\Delta^{1/\beta}W_{1}$ have the same law, it is intuitively
clear that $I_{\Delta}^{\sigma\sigma}(\sigma,\beta,0,\delta_{0})$ does not
depend on $\Delta$. Indeed, let us denote by $h_{\beta}$ the density of
$W_{1}$, which is defined through (\ref{eq:Wchar}). The density of $X_{\Delta
}=\sigma W_{\Delta}$ is
\[
p_{\Delta}(x|\sigma,\beta,0,\delta_{0})=\frac{1}{\sigma\Delta^{1/\beta}%
}~h_{\beta}\left(  \frac{x}{\sigma\Delta^{1/\beta}}\right)  .
\]
It is well known that $h_{\beta}$ is $C^{\infty}$ (by repeated integration of
the characteristic function), even, and that its $n-$th derivative $h_{\beta
}^{(n)}$ behaves as follows (the first two derivatives are denoted $h^{\prime
}$ and $h^{\prime\prime}$):
\begin{equation}
\left\vert h_{\beta}^{(n)}(w)\right\vert ~\sim~\left\{
\begin{array}
[c]{ll}%
{\frac{c_{\beta}(1+\beta)(2+\beta)\ldots(n-1+\beta)}{|w|^{n+1+\beta}}}\qquad &
\text{if \ }\beta<2\\[2.5mm]%
|w|^{n}e^{-w^{2}/2}/\sqrt{2\pi}\qquad & \text{if \ }\beta=2
\end{array}
\right.  \qquad\qquad\mbox{as }~~|w|\rightarrow\infty. \label{eq:dh^n}%
\end{equation}
where $c_{\beta}$ is given in (\ref{eq:cbeta}); this result follows from the
same series expansion as above). Let us also associate with $h_{\beta}$ the
following functions:
\begin{equation}
\breve{h}_{\beta}(w)=h_{\beta}(w)+wh_{\beta}^{\prime}(w),\qquad\widetilde
{h}_{\beta}(w)={\frac{\breve{h}_{\beta}(w)^{2}}{h_{\beta}(w)}}.
\label{eq:hfunctions}%
\end{equation}
Then $\tilde{h}_{\beta}$ is positive, even, continuous, and $\tilde{h}_{\beta
}(w)=0(1/|w|^{1+\beta})$ as $|w|\rightarrow\infty$, hence $\tilde{h}_{\beta}$
is Lebesgue--integrable.

Consider now
\begin{equation}
\mathcal{I}(\beta)=\int\widetilde{h}_{\beta}(w)dw, \label{eq:Icallig}%
\end{equation}
which is well defined and positive. Moreover if $\beta=2$ $(W$ is then
Brownian motion), $h_{2}$ is Gaussian and we have $h_{2}^{\prime}%
(w)=-wh_{2}(w)$, so $\widetilde{h}_{2}(w)= (1-w^{2}+w^{4})h_{2}(w)$ and
\begin{equation}
\beta=2~~~\Rightarrow~~~\mathcal{I}(\beta)=2. \label{eq:Icallig_beta=2}%
\end{equation}

The Fisher information for $\sigma$ associated with the observation of a
single variable $X_{\Delta}=\sigma W_{\Delta}$ is%
\[
I_{\Delta}^{\sigma\sigma}(\sigma,\beta,0,\delta_{0})=\int\frac{\left(
\partial_{\sigma}p_{\Delta}(x|\sigma,\beta,0,\delta_{0})\right)  ^{2}%
}{p_{\Delta}(x|\sigma,\beta,0,\delta_{0})}~dx=\frac{1}{\sigma^{3}%
\Delta^{1/\beta}}\int\widetilde{h}_{\beta}(x/\sigma\Delta^{1/\beta})~dx
\]
which, in light of (\ref{eq:Icallig}) and by a change of variable, reduces
to:
\begin{equation}
I_{\Delta}^{\sigma\sigma}(\sigma,\beta,0,\delta_{0})={\frac{1}{\sigma^{2}}%
}\mathcal{I}(\beta). \label{eq:Iss_Y=0}%
\end{equation}
So, as said before, this does not depend on $\Delta$. In fact, $\mathcal{I}%
(\beta)$ is simply the Fisher information at point $\sigma=1$ for the
statistical model in which we observe $\sigma W_{1}$ and $W_{1}$ is a variable
with density $h_{\beta}$.

\subsection{The index parameter $\beta.$}

Consider now the estimation of $\beta$. This problem was studied by
Dumouchel (1973a), who computed numerically the term $I_{\Delta}^{\beta
\beta}(\sigma,\beta,0,\delta_{0})$, including also an asymmetry parameter. It
is easily seen that $\beta\mapsto h_{\beta}(w)$ is differentiable on $(0,2]$,
and we denote by $\dot{h}_{\beta}(w)$ its derivative. However, instead of
(\ref{eq:dh^n}) one has
\begin{equation}
\left\vert \dot{h}_{\beta}(w)\right\vert ~\sim~\left\{
\begin{array}
[c]{ll}%
\frac{c_{\beta}\log|w|}{|w|^{1+\beta}}\qquad & \mbox{if
}~\beta<2\\[2.5mm]%
\frac{1}{|w|^{3}} & \mbox{if }~\beta=2,
\end{array}
\right.  \label{eq:h'beta}%
\end{equation}
as $|w|\rightarrow\infty$, by differentiation of the series expansion for the
stable density. Therefore the quantity
\begin{equation}
\mathcal{K}(\beta)=\int{\frac{\dot{h}_{\beta}(w)^{2}}{h_{\beta}(w)}}dw
\label{eq:Kcallig}%
\end{equation}
is finite when $\beta<2$ and infinite for $\beta=2$. This is the Fisher
information for estimating $\beta$, upon observing the single variable $W_{1}$.

Instead of computing the information quantities numerically based on
approximations of the stable density $h_{\beta}$, we study explicitly their
asymptotic behavior as $\Delta\rightarrow0.$ Excluding the degenerate case
where $\beta=2$ [see Dumouchel (1983) for the behavior of the MLE for
$\beta$ when $\beta=2$], the Fisher information for $\beta$ associated with
the observation of a single variable $\sigma W_{\Delta}$ when $\beta<2$ is%
\begin{align}
I_{\Delta}^{\beta\beta}(\sigma,\beta,0,\delta_{0})  &  =\int\frac{\left(
\partial_{\beta}p_{\Delta}(x|\sigma,\beta,0,\delta_{0})\right)  ^{2}%
}{p_{\Delta}(x|\sigma,\beta,0,\delta_{0})}~dx=\int\frac{\left(  \log
(\Delta)\breve{h}_{\beta}(w)+\beta^{2}\dot{h}_{\beta}(w)\right)  ^{2}}%
{\beta^{4}h_{\beta}(w)}~dw\nonumber\\
&  =\frac{(\log(\Delta))^{2}}{\beta^{4}}~\mathcal{I}(\beta)+\frac{2\log
(\Delta)}{\beta^{2}}\int\frac{\breve{h}_{\beta}(w)\dot{h}_{\beta}(w)}%
{h_{\beta}(w)}~dw+\mathcal{K}(\beta),\nonumber
\end{align}
and the middle integral in the last display is smaller than $\sqrt
{\mathcal{I}(\beta)\mathcal{K}(\beta)}$ by Cauchy--Schwarz. Therefore, as
$\Delta\rightarrow0,$ we have
\begin{equation}
{\frac{I_{\Delta}^{\beta\beta}(\sigma,\beta,0,\delta_{0})}{(\log
(1/\Delta))^{2}}}\rightarrow\frac{1}{\beta^{4}}~\mathcal{I}(\beta).
\label{eq:Ibb_Y=0_limit}%
\end{equation}

\subsection{The cross $(\sigma,\beta)\ $term.}

As for the cross-term, when of course $\beta<2$ again, we have%
\begin{align}
I_{\Delta}^{\sigma\beta}(\sigma,\beta,0,\delta_{0})  &  =\int\frac
{\partial_{\sigma}p_{\Delta}(x|\sigma,\beta,0,\delta_{0})~~\partial_{\beta
}p_{\Delta}(x|\sigma,\beta,0,\delta_{0})}{p_{\Delta}(x|\sigma,\beta
,0,\delta_{0})}~dx\nonumber\\
&  =-\frac{1}{\sigma\beta^{2}}\int\frac{\breve{h}_{\beta}(w)\left(
\log(\Delta)\breve{h}_{\beta}(w)-\beta^{2}\dot{h}_{\beta}(w)\right)
}{h_{\beta}(w)}~dw\nonumber\\
&  =-\frac{\log(\Delta)}{\sigma\beta^{2}}~\mathcal{I}(\beta)-\frac{1}{\sigma
}\int{\frac{\breve{h}_{\beta}(w)\dot{h}_{\beta}(w)}{h_{\beta}(w)}%
~dw.}\nonumber
\end{align}
Therefore, as $\Delta\rightarrow0,$ we have
\begin{equation}
{\frac{I_{\Delta}^{\sigma\beta}(\sigma,\beta,0,\delta_{0})}{\log(1/\Delta)}%
}\rightarrow\frac{1}{\sigma\beta^{2}}~\mathcal{I}(\beta).
\label{eq:Isb_Y=0_limit}%
\end{equation}

\subsection{The information for a translation model.}

We will see another information appear in some of the forthcoming formulas,
namely the Fisher information associated with the estimation of the real
number $a$ for the model where one observes the single variable $W_{1}+a$.
This Fisher information is of course the following number:
\begin{equation}
\mathcal{J}(\beta)=\int{\frac{h_{\beta}^{\prime}(w)^{2}}{h_{\beta}(w)}}dw.
\label{eq:Jcallig}%
\end{equation}
Observe in particular that
\begin{equation}
\beta=2~~~\Rightarrow~~~\mathcal{J}(\beta)=1. \label{eq:Jcallig_beta=2}%
\end{equation}

\subsection{Some consequences for the estimation.}

If now we come back to our setting where $n$ values (or increments) of
$X=\sigma W$ are observed along a time lag $\Delta_{n}$, we see that when
$\beta$ is known we can hope for estimators $\widehat{\sigma}_{n}$ for
$\sigma$ which are asymptotically efficient in the sense that $\sqrt{n}~
(\widehat{\sigma}_{n}-\sigma)$ converges in law to ${\mathcal{N}}0,\sigma
^{2}/{\mathcal{I}}(\beta))$, whatever $\Delta_{n}$ behaves like as $n\to
\infty$, and of course the MLE satisfies that.

When it comes to estimating $\beta$, things are different. When $\Delta
_{n}\rightarrow0$ and when the true value is $\beta<2$, we can hope for
estimators converging to $\beta$ with the faster rate $\sqrt{n}~\log
(1/\Delta_{n})$ (and we provide such estimators in A\"\i t--Sahalia
and Jacod~(2004)). Some
of these estimators are constructed using the jumps of the process of a size
greater than some threshold, as in H\"opfner and Jacod (1994). Note
also that if we
suspect that $\beta=2$ we would rather perform a test, as advised by
Dumouchel (1973a), and anyway in this case the behavior of the Fisher
information does not provide much insight.

\section{The general semiparametric case.\label{sec:semiparametric}}

The data generating process is now given by (\ref{eq:X=sW+aY}). We are
interested in estimating $(\sigma,\beta)$, and in some instances $\theta$ as
well, leaving the distribution $G\in\mathcal{G}_{\beta}$ unspecified.

\subsection{Estimation of $(\sigma,\beta)$.}

We start by studying whether the limiting behavior of $I_{\Delta}%
^{\sigma\sigma}(\sigma,\beta,\theta,G)$ when $\beta=2$ and of the
$(\sigma,\beta)$ block of the matrix $I_{\Delta}(\sigma,\beta,\theta,G)$ when
$\beta<2$ is affected by the presence of $Y.$ First, we have the intuitively
obvious majoration of Fisher's information in presence of $Y$ by the one for
which $Y$ is absent. Note that in this result no assumption whatsoever is made
on $Y$ (except of course that it is independent of $W$):

\begin{theorem}
\label{theorem:generalmajoration} For any $\Delta>0$ we have
\begin{equation}
I_{\Delta}^{\sigma\sigma}(\sigma,2,\theta,G)\leq I_{\Delta}^{\sigma\sigma
}(\sigma,2,0,\delta_{0}) \label{FI-MAJ}%
\end{equation}
and, when $\beta<2$, the difference
\[
\left(
\begin{array}
[c]{ccc}%
I_{\Delta}^{\sigma\sigma}(\sigma,\beta,0,\delta_{0}) & ~~ & I_{\Delta}%
^{\sigma\beta}(\sigma,\beta,0,\delta_{0})\\
I_{\Delta}^{\sigma\beta}(\sigma,\beta,0,\delta_{0}) &  & I_{\Delta}%
^{\beta\beta}(\sigma,\beta,0,\delta_{0})
\end{array}
\right)  -\left(
\begin{array}
[c]{ccc}%
I_{\Delta}^{\sigma\sigma}(\sigma,\beta,\theta,G) & ~~ & I_{\Delta}%
^{\sigma\beta}(\sigma,\beta,\theta,G)\\
I_{\Delta}^{\sigma\beta}(\sigma,\beta,\theta,G) &  & I_{\Delta}^{\beta\beta
}(\sigma,\beta,\theta,G)
\end{array}
\right)
\]
is a positive semi--definite matrix, and in particular we have:
\begin{equation}
I_{\Delta}^{\beta\beta}(\sigma,\beta,\theta,G)\leq I_{\Delta}^{\beta\beta
}(\sigma,\beta,0,\delta_{0}). \label{FI-MAJB}%
\end{equation}

\end{theorem}

Next, how does the limit as $\Delta\rightarrow0$ of $I_{\Delta}(\sigma
,\beta,\theta,G)$ compare to that of $I_{\Delta}(\sigma,\beta,0,\delta_{0})$?
For instance, given that in the absence of $Y$ we can estimate $\sigma$ with
information $\mathcal{I}(\beta)/\sigma^{2},$ we would like to find out what is
the impact, if any, of the presence of $Y$ on the information we can gather
about that parameter from the discrete observations where $W$ is perturbed by
$Y:$
\[
\chi_{i}^{n}=\sigma(W_{i\Delta_{n}}-W_{(i-1)\Delta_{n}})+\theta(Y_{i\Delta
_{n}}-Y_{(i-1)\Delta_{n}})
\]
for $i=1,\ldots,n.$

The answer to that question is given by the following.

\begin{theorem}
\label{theorem:sigma&beta} a) If $G\in\mathcal{G}_{\beta}$ we have as
$\Delta\rightarrow0$:
\begin{equation}
I_{\Delta}^{\sigma\sigma}(\sigma,\beta,\theta,G)\rightarrow\frac{1}{\sigma
^{2}}\mathcal{I}(\beta), \label{FI-IN}%
\end{equation}
and also, when $\beta<2$:
\begin{equation}
{\frac{I_{\Delta}^{\beta\beta}(\sigma,\beta,\theta,G)}{(\log(1/\Delta)^{2}}%
}\rightarrow\frac{1}{\beta^{4} }\mathcal{I}(\beta),\qquad{\frac{I_{\Delta
}^{\sigma\beta}(\sigma,\beta,\theta,G)}{\log(1/\Delta)}}\rightarrow\frac
{1}{\sigma\beta^{2}} \mathcal{I}(\beta) \label{FI-IN'}%
\end{equation}

b) For any $\phi\in\Phi$ and $\alpha\in(0,\beta\rbrack$ and $K>0$, we have as
$\Delta\rightarrow0$:
\begin{equation}
\left\{
\begin{array}
[c]{l}%
\sup_{G\in\mathcal{G}(\phi,\alpha),|\theta|\leq K}\left|  I_{\Delta}%
^{\sigma\sigma}(\sigma,\beta,\theta,G) -{\frac{\mathcal{I}(\beta)}{\sigma^{2}%
}}\right|  \rightarrow0,\\[2.5mm]%
\beta<2\quad\Rightarrow\quad\left\{
\begin{array}
[c]{l}%
\sup_{G\in\mathcal{G}(\phi,\alpha),|\theta|\leq K}\left|  \frac{I_{\Delta
}^{\beta\beta}(\sigma,\beta,\theta,G)}{(\log(1/\Delta))^{2}} -{\frac
{\mathcal{I}(\beta)}{\beta^{4}}}\right|  \rightarrow0,\\[2.5mm]%
\sup_{G\in\mathcal{G}(\phi,\alpha),|\theta|\leq K}\left|  \frac{I_{\Delta
}^{\sigma\sigma}(\sigma,\beta,\theta,G)}{\log(1/\Delta)} -{\frac
{\mathcal{I}(\beta)}{\sigma\beta^{2}}}\right|  \rightarrow0.
\end{array}
\right.
\end{array}
\right.  \label{eq:Iss_uniform}%
\end{equation}

c) For each $n,$ let $G^{n}$ be the standard symmetric stable law of index
$\alpha_{n}$, with $\alpha_{n}$ a sequence strictly increasing to $\beta$.
Then for any sequence $\Delta_{n}\rightarrow0$ such that $(\beta-\alpha
_{n})\log\Delta_{n}\rightarrow0$ (i.e. the rate at which $\Delta
_{n}\rightarrow0$ is slow enough), the sequence of numbers $I_{\Delta_{n}%
}^{\sigma\sigma}(\sigma,\beta,\theta,G^{n})$ (resp. $I_{\Delta_{n}}%
^{\beta\beta}(\sigma,\beta,\theta,G^{n})/(\log(1/\Delta_{n}))^{2}$ when
further $\beta<2$) converges to a limit which is strictly less than
$\mathcal{I}(\beta)/\sigma^{2}$ (resp. $\mathcal{I}(\beta)/\beta^{4}$).
\end{theorem}

In other words, at their respective leading orders in $\Delta,$ the presence
of $Y\ $has no impact on the information terms $I_{\Delta}^{\sigma\sigma},$
$I_{\Delta}^{\beta\beta}$ and $I_{\Delta}^{\sigma\beta}$, as soon as $Y$ is
\textquotedblleft dominated\textquotedblright\ by $W$: so, in the limit where
$\Delta\rightarrow0,$ the parameters $\sigma$ and $\beta$ can be estimated
with the exact same degree of precision whether $Y\ $is present or not.
Moreover, part (b) states the convergence of Fisher's information is
\emph{uniform} on the set $\mathcal{G}(\phi,\alpha)$ and $|\theta|\leq K$ for
all $\alpha\in\lbrack0,\beta]$; this settles the case where $G$ and $\theta$
are considered as nuisance parameters when we estimate $\sigma$ and $\beta$.

But as $\alpha$ tends to $\beta,$ the \emph{convergence disappears}, as stated
in part (c). This shows that the class $\mathcal{G}_{\beta}$ is effectively
the largest one for which the presence of a $Y$ process does not affect the
estimation of the parameters of the process $\sigma W.$ For example, if
$\beta=2,$ in part (c)\ take $G^{n}$ to be the symmetric stable law with index
$\alpha_{n}\in(0,2)$ and scale parameter $s$ in the sense that its
characteristic function is $u\mapsto\exp\left(  -\frac{s^{2}}{2}%
|u|^{\alpha_{n}}\right)  $. Then if $\alpha_{n}\rightarrow2$, for all
sequences $\Delta_{n}\rightarrow0$ satisfying $(2-\alpha_{n})\log\Delta
_{n}\rightarrow0$ we have
\[
I_{\Delta_{n}}^{\sigma\sigma}(\sigma,\beta,1,G^{n})\rightarrow{\frac{2}%
{\sigma^{2}+s^{2}}}.
\]
This is of course to be expected, since in the limit we are observing
$\sqrt{\sigma^{2}+s^{2}}~W$, and we supposedly know $s$ and wish to estimate
$\sigma$.

Another interesting feature, due to the fact that $\Delta\rightarrow0$, is
that the limiting behavior of $I_{\Delta}^{\beta\beta}$ and $I_{\Delta
}^{\sigma\beta}$ when $\beta<2$ involves $\mathcal{I}(\beta)$ but not
$\mathcal{K}(\beta)$, as one could have guessed at first glance.

\subsection{Estimation of $\theta$.}

For the entries of Fisher's information matrix involving the parameter
$\theta$, things are more complicated. First, observe that $I_{\Delta}%
^{\theta\theta}(0,\beta,\theta,G)$ (that is the Fisher information for the
model $X=\theta Y$) does not necessarily exist, but of course if it does we
have an inequality similar to (\ref{FI-MAJ}) for all $\sigma$:
\begin{equation}
\label{FI-MAJ-Y}I_{\Delta}^{\theta\theta}(\sigma,\beta,\theta,G)\leq
I_{\Delta}^{\theta\theta}(0,\beta,\theta,G).
\end{equation}

Contrary to (\ref{FI-MAJ}), however, this is a very rough estimate, which does
not take into account the properties of $W$. The $(\theta,\theta)$--Fisher
information is usually much smaller than what the right side above suggests,
and we give below a more accurate estimate when $Y$ has second moments, but
without the \textquotedblleft domination\textquotedblright\ assumption that
$G\in\mathcal{G}_{\beta}$. Recall the notation (\ref{eq:Jcallig}).

\begin{theorem}
\label{theorem:theta} If $Y_{1}$ has a finite variance $v$ and a mean $m$, we
have
\begin{equation}
I_{\Delta}^{\theta\theta}(\sigma,\beta,\theta,G)\leq\frac{\mathcal{J}(\beta
)}{\sigma^{2}}~\left(  m^{2}\Delta^{2-2/\beta}+v\Delta^{1-2/\beta}\right)  .
\label{FI-INtheta1}%
\end{equation}

\end{theorem}

This estimate holds for all $\Delta>0$. The asymptotic variant, which says
that
\begin{equation}
\label{FI-INtheta2}\limsup_{\Delta\to0}~\Delta^{2/\beta-1}~ I_{\Delta}%
^{\theta\theta}(\sigma,\beta,\theta,G)\leq\frac{\delta~\mathcal{J}(\beta
)}{\sigma^{2}},
\end{equation}
is sharp in some cases and not in others, as we will see in the examples
below. These examples will also enlighten the fact that the ``translation''
Fisher information $\mathcal{J}(\beta)$ comes into the picture here.

\section{Examples.\label{sec:examples}}

The calculations of the previous section involving the parameter $\theta$ can
be made fully explicit if we specify the distribution of the process $Y$, in
some cases at least. We will always suppose that $\beta$ is known in these examples.

\subsection{Stable process plus drift.\label{sec:drift}}

Here we assume that $Y_{t}=t$, so $G_{\Delta}=\delta_{\Delta}$ and
$G=\delta_{1}$ (recall that the notation $\delta$ means a Dirac mass):

\begin{theorem}
\label{theorem:S+D} The $2\times2$ Fisher information matrix for estimating
$(\sigma,\theta)$ is
\begin{equation}
\left(
\begin{array}
[c]{ccc}%
I_{\Delta}^{\sigma\sigma}(\sigma,\beta,\theta,\delta_{1}) & ~~ & I_{\Delta
}^{\sigma\theta}(\sigma,\beta,\theta,\delta_{1})\\[2mm]%
I_{\Delta}^{\sigma a}(\sigma,\beta,\theta,\delta_{1}) &  & I_{\Delta}%
^{aa}(\sigma,\beta,\theta,\delta_{1})
\end{array}
\right)  ={\frac{1}{\sigma^{2}}}~\left(
\begin{array}
[c]{cc}%
\mathcal{I}(\beta)\quad & 0\\[2mm]%
0\quad & \Delta^{2-2/\beta}~\mathcal{J}(\beta)
\end{array}
\right)  . \label{eq:Fisher_S+D}%
\end{equation}

\end{theorem}

This has several interesting consequences (we will denote by $T_{n}%
=n\Delta_{n}$ the length of the \emph{observation window}):

\begin{enumerate}
\item If $\theta$ is known, one may hope for estimators $\widehat{\sigma}_{n}$
for $\sigma$ such that $\sqrt{n}(\widehat{\sigma}_{n}-\sigma)~\overset
{d}{\longrightarrow}~\mathcal{N}(0,\sigma^{2}/\mathcal{I}(\beta))$ (that is,
\emph{asymptotically efficient}\ in the Cramer--Rao sense). As a matter of
fact, in this setting, observing $\chi_{i}^{n}$ is equivalent to observing
$\chi^{\prime}{}_{i}^{n}=\chi_{i}^{n}-a\Delta_{n}$, so we are in the situation
of Section \ref{sec:base}.

\item If $\sigma$ is known, one may hope for estimators $\widehat{\theta}_{n}$
for $\theta$ such that $\sqrt{n}\Delta_{n}^{1-1/\beta}(\widehat{\theta}%
_{n}-\theta)$ converges in law to $\mathcal{N}(0,\sigma^{2}/\mathcal{J}%
(\beta))$. If $\beta=2$ the rate is thus $\sqrt{T_{n}}$: this is in accordance
with the well known fact that for a diffusion the rate for estimating the
drift coefficient is the square root of the total observation window, that is
$\sqrt{T_{n}}$ here; moreover in this case, the variable $X_{T_{n}}/T_{n}$ is
$\mathcal{N}(\theta,\sigma^{2}/T_{n})$; so $\widehat{\theta}_{n}=X_{T_{n}%
}/T_{n}$ is an asymptotically efficient estimator for $\theta$ (recall that
$\mathcal{J}(\beta)=1$ when $\beta=2$). When $\beta<2$ we have $1-1/\beta
<1/2$, so the rate is bigger than $\sqrt{T_{n}}$, and it increases when
$\beta$ decreases; when $\beta<1$ this rate is even bigger than $\sqrt{n}$.

Observe that here $Y_{1}$ has mean $m=1$ and variance $v=0$: so the estimate
(\ref{FI-INtheta1}) is indeed an equality. The fact that the translation
Fisher information $\mathcal{J}(\beta)$ appears here is transparent.

\item If both $\sigma$ and $\theta$ are unknown, one may hope for estimators
$\widehat{\sigma}_{n}$ and $\widehat{\theta}_{n}$ such that the pairs
$(\sqrt{n}(\widehat{\sigma}_{n}-\sigma),\sqrt{n}\Delta_{n}^{1-1/\beta
}(\widehat{\theta}_{n}-\theta))$ converge in law to the product $\mathcal{N}%
(0,\sigma^{2}/\mathcal{I}(\beta))\otimes\mathcal{N}(0,\sigma^{2}%
/\mathcal{J}(\beta))$.
\end{enumerate}

\subsection{Stable process plus Poisson process.\label{sec:Poisson}}

Here we assume that $Y$ is a standard Poisson process\ (jumps of size $1$,
intensity $1$), whose law we write as $G=P$. We can describe the limiting
behavior of the $(\sigma,\theta)$ block of the matrix $I_{\Delta}(\sigma
,\beta,\theta,P)$ as $\Delta\rightarrow0$.

\begin{theorem}
\label{theorem:S+P} If $Y$ is a standard Poisson process we have, as
$\Delta\rightarrow0$:
\begin{align}
I_{\Delta}^{\sigma\sigma}(\sigma,\beta,\theta,P)~  &  \rightarrow~{\frac
{1}{\sigma^{2}}}~\mathcal{I}(\beta)\label{eq:Iss_Y=Poisson}\\
\Delta^{1/\beta-1/2}~I_{\Delta}^{\sigma\theta}(\sigma,\beta,\theta,P)~  &
\rightarrow~0\label{eq:Isa_Y=Poisson}\\
\Delta^{2/\beta-1}~I_{\Delta}^{\theta\theta}(\sigma,\beta,\theta,P)~  &
\rightarrow~{\frac{1}{\sigma^{2}}}~\mathcal{J}(\beta) \label{eq:Iaa_Y=Poisson}%
\end{align}

\end{theorem}

Since $P\in\mathcal{G}_{\beta}$, (\ref{eq:Iss_Y=Poisson}) is nothing else than
the first part of (\ref{FI-IN}). One could prove more than
(\ref{eq:Isa_Y=Poisson}), namely that $\sup_{\Delta}\Delta^{1/\beta-1}~
|I_{\Delta}^{\sigma\theta}(\sigma,\beta,\theta,P)|\leq\infty$. Here again, we
deduce some interesting consequences:

\begin{enumerate}
\item If $\sigma$ is known, one may hope for estimators $\widehat{\theta}_{n}$
for $\theta$ such that $\sqrt{n\Delta_{n}^{1-2/\beta}}~(\widehat{\theta}%
_{n}-\theta)$ converge in law to $\mathcal{N}(0,\sigma^{2}/\mathcal{J}%
(\beta))$. So the rate is \emph{bigger}\ than $\sqrt{n}$, except when
$\beta=2$. More generally, if both $\sigma$ and $\theta$ are unknown, one may
hope for estimators $\widehat{\sigma}_{n}$ and $\widehat{\theta}_{n}$ such
that the pairs $\left(  \sqrt{n}~(\widehat{\sigma}_{n}-\sigma),\sqrt
{n\Delta_{n}^{1-2/\beta}}~ (\widehat{\theta}_{n}-\theta)\right)  $ converge in
law to the product $\mathcal{N}(0,\sigma^{2}/\mathcal{I}(\beta))\otimes
\mathcal{N}(0,\sigma^{2} /\mathcal{J}(\beta))$.

\item However, the above--described behavior of any estimator $\widehat
{\theta}_{n}$ cannot be true when $T_{n}=n\Delta_{n}$ does not go to infinity,
because in this case there is a positive probability that $Y$ has no jump on
the biggest observed interval, and so no information about $\theta$ can be
drawn from the observations in that case. It is true, though, when
$T_{n}\rightarrow\infty$, because $Y$ will eventually have infinitely many
jumps on the observed intervals. This discrepancy between the asymptotic
behavior of Fisher information and of estimators shows that some care must be
taken when the Fisher information is used as a measure of the quality of estimators.

\item Observe that here $Y_{1}$ has mean $m=1$ and variance $v=1$. So in view
of (\ref{eq:Iaa_Y=Poisson}) the asymptotic estimate (\ref{FI-INtheta2}) is sharp.
\end{enumerate}

\subsection{Stable process plus compound Poisson
process.\label{sec:compoundPoisson}}

Here we assume that $Y$ is a compound Poisson process with arrival rate
$\lambda$ and law of jumps $\mu$: that is, the characteristics of $G$ are
$b=\lambda\int_{\{|x|\leq1\}}x\mu(dx)$ and $c=0$ and $F=\lambda\mu$. We then
write $G=P_{\lambda,\mu}$, which belongs to $\mathcal{G}_{\beta}$.

We will further assume that $\mu$ has a density $f$ satisfying:
\begin{equation}
\lim_{|u|\rightarrow\infty}uf(u)=0,\qquad\sup_{u}(|f^{\prime}%
(u)|(1+|u|))<\infty. \label{eq:cond-f}%
\end{equation}
We also suppose that the \textquotedblleft multiplicative\textquotedblright%
\ Fisher information associated with $\mu$ (that is, the Fisher information
for estimating $\theta$ in the model when one observes a single variable
$\theta U$ with $U$ distributed according to $\mu$) exists. It then has the
form
\begin{equation}
\mathcal{L}=\int\frac{(uf^{\prime}(u)+f(u))^{2}}{f(u)}~du. \label{eq:FishMult}%
\end{equation}

We can describe the limiting behavior of the $(\sigma,\theta)$ block of the
matrix $I_{\Delta}(\sigma,\beta,\theta,P_{\lambda,\mu})$ as $\Delta
\rightarrow0$.

\begin{theorem}
\label{theorem:S+CP} If $Y$ is a compound Poisson process satisfying
(\ref{eq:cond-f}) and such that $\mathcal{L}$ in (\ref{eq:FishMult}) is
finite, we have as $\Delta\rightarrow0$:
\begin{align}
I_{\Delta}^{\sigma\sigma}(\sigma,\beta,\theta,P_{\lambda,\mu})~  &
\rightarrow~{\frac{1}{\sigma^{2}}}~\mathcal{I}(\beta)
\label{eq:Iss_Y=CPoisson}\\
~\frac{1}{\sqrt{\Delta}}~I_{\Delta}^{\sigma\theta}(\sigma,\beta,\theta
,P_{\lambda,\mu})~  &  \rightarrow~0 \label{eq:Isa_Y=CPoisson}%
\end{align}
and
\begin{equation}
\frac{\lambda^{2}}{\theta^{2}}~\int\frac{(xf^{\prime}(x)+f(x))^{2}}{\lambda
f(x)+\frac{c_{\beta}~\sigma^{\beta}}{\theta^{\beta}|x|^{1+\beta}}}%
~dx\leq\liminf\frac{1}{\Delta}~I_{\Delta}^{\theta\theta}(\sigma,\beta
,\theta,P_{\lambda,\mu})\leq\limsup\frac{1}{\Delta}~I_{\Delta}^{\theta\theta
}(\sigma,\beta,\theta,G)~\leq{\frac{1}{\theta^{2}}}~\mathcal{L}
\label{eq:Iaa_Y=CPoisson}%
\end{equation}
when $\beta<2$ ($c_{\beta}$ is the constant defined in (\ref{eq:cbeta})), and
also, when $\beta=2$:
\begin{equation}
\frac{1}{\Delta}~I_{\Delta}^{\theta\theta}(\sigma,\beta,\theta,P_{\lambda,\mu
})\rightarrow{\frac{1}{\theta^{2}}}~\mathcal{L}. \label{eq:Iaa2_Y=CPoisson}%
\end{equation}

\end{theorem}

As for the previous theorem, (\ref{eq:Iss_Y=CPoisson}) is nothing else than
the first part of (\ref{FI-IN}). We could prove more than
(\ref{eq:Isa_Y=CPoisson}), namely that $\sup_{\Delta}\frac{1}{\Delta
}~|I_{\Delta}^{\sigma\theta}(\sigma,\beta,\theta,P_{\lambda,\mu})|<\infty$.
Here again, we deduce some interesting consequences:

\begin{enumerate}
\item One may hope for estimators $\widehat{\theta}_{n}$ for $\theta$ such
that $\sqrt{T_{n}}(\widehat{\theta}_{n}-\theta)$ is tight (the rate is the
same as for the case $Y_{t}=t$), and is even asymptotically normal when
$\beta=2$.

\item However, this is not true when $T_{n}$ does not go to infinity, for the
same reason as for the previous theorem.

\item When the measure $\mu$ has a second order moment, the right side of
(\ref{FI-INtheta1}) is larger than the result of the previous theorem, so the
estimate in Theorem \ref{theorem:theta} is not sharp.
\end{enumerate}

The rates for estimating $\theta$ in the two previous theorems, and the
limiting Fisher information as well, can be explained as follows (supposing
that $\sigma$ is known and that we have $n$ observations and that
$T_{n}\rightarrow\infty$):

\begin{enumerate}
\item For Theorem \ref{theorem:S+P}: $\theta$ comes into the picture whenever
the Poisson process has a jump. On the interval $[0,T_{n}]$ we have an average
of $T_{n}$ jumps, most of them being isolated in an interval $(i\Delta
_{n},(i+1)\Delta_{n}]$. So it essentially amounts to observing $T_{n}$ (or
rather the integer part $[T_{n}]$) independent variables, all distributed as
$\sigma\Delta_{n}^{1/\beta}W_{1}+\theta$. The Fisher information for each of
those (for estimating $\theta$) is $J(\beta)/\sigma^{2}\Delta^{2/\beta}$, and
the \textquotedblleft global\textquotedblright\ Fisher information, namely
$nI_{\Delta_{n}}^{\theta\theta}$, is approximately $T_{n}J(\beta)/\sigma
^{2}\Delta_{n}^{2/\beta}\sim J(\beta)/\sigma^{2}\Delta_{n}^{2/\beta-1}$.

\item For Theorem \ref{theorem:S+CP}: Again $\theta$ comes into the picture
whenever the compound Poisson process has a jump. We have an average of
$\lambda T_{n}$ jumps, so it essentially amounts to observing $\lambda T_{n}$
independent variables, all distributed as $\sigma\Delta_{n}^{1/\beta}%
W_{1}+\theta V$ where $V$ has the distribution $\mu$. The Fisher information
for each of those (for estimating $\theta$) is approximately $L/\theta^{2}$
(because the variable $\sigma\Delta_{n}^{1/\beta}W_{1}$ is negligible), and
the \textquotedblleft global\textquotedblright\ Fisher information
$nI_{\Delta_{n}}^{\theta\theta}$ is approximately $\lambda T_{n}L\theta
^{2}\sim n\Delta_{n}L/\theta^{2}$. This explains the rate in
(\ref{eq:Iaa_Y=CPoisson}), and is an indication that (\ref{eq:Iaa2_Y=CPoisson}%
) may be true even when $\beta<2$, although we have been unable to prove it
thus far. \bigskip
\end{enumerate}

\subsection{Two stable processes.\label{sec:twostables}}

Our last example is about the case where $Y$ is also a symmetric stable
process with index $\alpha,$ $\alpha<\beta$. We write $G=S_{\alpha}.$
Surprisingly, the results are quite involved, in the sense that for estimating
$\theta$ we have different situations according to the relative values of
$\alpha$ and $\beta$. We obviously still have (\ref{eq:Iss_Y=CPoisson}), so we
concentrate on the term $I_{\Delta}^{\theta\theta}$ and ignore the cross term
in the statement of the following theorem:

\begin{theorem}
\label{theorem:S+S} If $Y$ is a standard symmetric stable process with index
$\alpha<\beta$, we have as $\Delta\rightarrow0$:
\begin{align}
\beta=2\quad &  \Rightarrow\quad\frac{(\log(1/\Delta))^{\alpha/2}}%
{\Delta^{\frac{\beta-\alpha}{\beta}}}~I_{\Delta}^{\theta\theta}(\sigma
,\beta,\theta,S_{\alpha})~\rightarrow~\frac{2\alpha c_{\alpha}\beta^{\alpha
/2}}{\theta^{2-\alpha}\sigma^{\alpha}(2(\beta-\alpha))^{\alpha/2}%
}\label{eq:SS1}\\
\beta<2,~\alpha>\frac{\beta}{2}\quad &  \Rightarrow\quad\frac{1}{\Delta
^{\frac{2(\beta-\alpha)}{\beta}}}~I_{\Delta}^{\theta\theta}(\sigma
,\beta,\theta,S_{\alpha})~\nonumber\\
&  \rightarrow\frac{\alpha^{2}c_{\alpha}^{2}\theta^{2\alpha-2}}{\sigma
^{2\alpha}}\int\frac{\left(  \int_{\mathbf{R}}|y|^{1-\alpha}dy\int_{0}%
^{1}(1-v)h_{\beta}^{\prime\prime}(x-yv)dv\right)  ^{2}}{h_{\beta}%
(y)}~dx\label{eq:SS2}\\
\beta<2,~\alpha=\frac{\beta}{2}\quad &  \Rightarrow\quad\frac{1}{\Delta
^{\frac{2(\beta-\alpha)}{\beta}}\log(1/\Delta)}~I_{\Delta}^{\theta\theta
}(\sigma,\beta,\theta,S_{\alpha})~\rightarrow~\frac{2\alpha(\beta
-\alpha)c_{\alpha}^{2}\theta^{2\alpha-2}}{\beta c_{\beta}\sigma^{2\alpha}%
}\label{eq:SS3}\\
\beta<2,~\alpha<\frac{\beta}{2}\quad &  \Rightarrow\quad\frac{1}{\Delta
}~I_{\Delta}^{\theta\theta}(\sigma,\beta,\theta,S_{\alpha})~\rightarrow
~\frac{\alpha^{2}c_{\alpha}^{2}\theta^{2\alpha-2}}{\sigma^{2\alpha}}~\int
\frac{1}{c_{\beta}|z|^{1+2\alpha-\beta}+c_{\alpha}\theta^{\alpha}%
|z|^{1+\alpha}/\sigma^{\alpha}}~dz \label{eq:SS4}%
\end{align}

\end{theorem}

Then if $\sigma$ is known one may hope to find estimators $\widehat{\theta
}_{n}$ for $\theta$ such that $u_{n}(\widehat{\theta}_{n}-\theta)$ converges
in law to $N(0,V)$, with%

\[%
\begin{array}
[c]{ll}%
u_{n}=\frac{\sqrt{n}~\Delta_{n}^{\frac{\beta-\alpha}{2\beta}}}{(\log
(1/\Delta_{n})^{\alpha/4}}\qquad & \text{if \ }~\beta=2\\[2.5mm]%
u_{n}=\sqrt{n}~\Delta_{n}^{\frac{\beta-\alpha}{\beta}}\qquad & \text{if
\ }~\beta<2,~\alpha>\beta/2\\[2.5mm]%
u_{n}=\sqrt{n}~\Delta_{n}^{\frac{\beta-\alpha}{\beta}}\sqrt{\log(1/\Delta
_{n})}\qquad & \text{if \ }~\beta<2,~\alpha=\beta/2\\[2.5mm]%
u_{n}=\sqrt{T_{n}}\qquad & \text{if \ }~\beta<2,~\alpha<\beta/2.
\end{array}
\]
and of course the asymptotic variance $V$ should be the inverse of the right
hand sides in (\ref{eq:SS1})--(\ref{eq:SS4}).

\section{The multiplicative model.\label{sec:multiplicative}}

Another interesting situation is the model (\ref{eq:X=sW+sY}), when $\beta$ is
known. If we observe a single variable $X_{\Delta}$, the corresponding Fisher
information is obviously
\[
I_{\Delta}^{\prime}(\sigma,\beta,G)=I_{\Delta}^{\sigma\sigma}(\sigma
,\beta,\sigma,G)+2I_{\Delta}^{\sigma\theta}(\sigma,\beta,\sigma,G)+I_{\Delta
}^{\theta\theta}(\sigma,\beta,\sigma,G).
\]
We will not develop a full theory here, but we translate the examples of the
previous section in this setting. In view of the previous results, the proofs
of the next three theorem follow, and these results show the variety of
situations we may encounter for this multiplicative model. The problems dealt
with in Theorems \ref{theorem:S+P-M}-\ref{theorem:S+CP-M} below have been
solved previously by Far (2001) and Jedidi (2001).

\begin{theorem}
\label{theorem:S+D-M} If $Y_{t}=t$ the Fisher information for estimating
$\sigma$ in the model (\ref{eq:X=sW+sY}) satisfies, as $\Delta\rightarrow0$:
\begin{align}
I_{\Delta}^{\prime}(\sigma,\beta,\delta_{1})  &  \rightarrow\frac{1}%
{\sigma^{2}}~\mathcal{I}(\beta) & \mbox{if }~\beta &  \in(1,2]\nonumber\\
I_{\Delta}^{\prime}(\sigma,\beta,\delta_{1})  &  \rightarrow\frac{1}%
{\sigma^{2}}~(\mathcal{I}(\beta)+\mathcal{J}(\beta)) & \mbox{if }~\beta &
=1\nonumber\\
\Delta^{2/\beta-2}I_{\Delta}^{\prime}(\sigma,\beta,\delta_{1})  &
\rightarrow\frac{1}{\sigma^{2}}~\mathcal{J}(\beta) & \mbox{if }~\beta &
\in(0,1).\nonumber
\end{align}

\end{theorem}

So in this situation we may hope for estimators $\widehat{\sigma}_{n}$ for
$\sigma$ such that
\[%
\begin{array}
[c]{ll}%
\sqrt{n}(\widehat{\sigma}_{n}-\sigma)~\overset{d}{\longrightarrow}%
~\mathcal{N}\left(  0,{\frac{\sigma^{2}}{\mathcal{I}(\beta)}}\right)  \qquad &
\text{if \ }~\beta\in(1,2]\\[2.5mm]%
\sqrt{n}(\widehat{\sigma}_{n}-\sigma)~\overset{d}{\longrightarrow}%
~\mathcal{N}\left(  0,{\frac{\sigma^{2}}{\mathcal{I}(\beta)+\mathcal{J}%
(\beta)}}\right)  \qquad & \text{if \ }~\beta=1\\[2.5mm]%
\sqrt{n}\Delta_{n}^{1-1/\beta}(\widehat{\sigma}_{n}-\sigma)~\overset
{d}{\longrightarrow}~\mathcal{N}\left(  0,{\frac{\sigma^{2}}{\mathcal{J}%
(\beta)}}\right)  \qquad & \text{if \ }~\beta\in(0,1).
\end{array}
\]

\begin{theorem}
\label{theorem:S+P-M} If $Y$ is a standard Poisson process, the Fisher
information for estimating $\sigma$ in the model (\ref{eq:X=sW+sY}) satisfies,
as $\Delta\rightarrow0$:
\begin{align}
I_{\Delta}^{\prime}(\sigma,\beta,P)  &  \rightarrow\frac{1}{\sigma^{2}%
}~(\mathcal{I}(\beta)+\mathcal{J}(\beta)) & \mbox{if }~\beta &  =2\nonumber\\
\Delta^{2/\beta-1}I_{\Delta}^{\prime}(\sigma,\beta,P)  &  \rightarrow\frac
{1}{\sigma^{2}}~\mathcal{J}(\beta) & \mbox{if }~\beta &  \in(0,2).\nonumber
\end{align}

\end{theorem}

So in this situation, and as soon as $T_{n}\rightarrow\infty$, we may hope for
estimators $\widehat{\sigma}_{n}$ for $\sigma$ such that
\[%
\begin{array}
[c]{ll}%
\sqrt{n}(\widehat{\sigma}_{n}-\sigma)~\overset{d}{\longrightarrow}%
~\mathcal{N}\left(  0,{\frac{\sigma^{2}}{\mathcal{I}(\beta)+\mathcal{J}%
(\beta)}}\right)  \qquad & \text{if \ }~\beta=2\\[2.5mm]%
\sqrt{n\Delta_{n}^{1-2/\beta}}~(\widehat{\sigma}_{n}-\sigma)~\overset
{d}{\longrightarrow}~\mathcal{N}\left(  0,{\frac{\sigma^{2}}{\mathcal{J}%
(\beta)}}\right)  \qquad & \text{if \ }~\beta\in(0,2).
\end{array}
\]

\begin{theorem}
\label{theorem:S+CP-M} If $Y$ is a compound Poisson process satisfying
(\ref{eq:cond-f}) and such that $L$ in (\ref{eq:FishMult}) is finite, and also
if $Y$ is a symmetric stable process with index $\alpha<\beta$, the Fisher
information for estimating $\sigma$ in the model (\ref{eq:X=sW+sY}) satisfies,
as $\Delta\rightarrow0$:
\[
I_{\Delta}^{\prime}(\sigma,\beta,G)\rightarrow\frac{1}{\sigma^{2}}%
~\mathcal{I}(\beta).
\]

\end{theorem}

So in this situation, and as soon as $T_{n}\rightarrow\infty$, we may hope for
estimators $\widehat{\sigma}_{n}$ for $\sigma$ such that
\[
\sqrt{n}(\widehat{\sigma}_{n}-\sigma)~\overset{d}{\longrightarrow}%
~\mathcal{N}\left(  0,{\frac{\sigma^{2}}{\mathcal{J}(\beta)}}\right)  .
\]

\section{Proofs.\label{sec:proofs}}

\subsection{Preliminaries about the class $G(\phi,\alpha)$.
\label{ssec:technical}}

In the sequel, we denote by $C_{\gamma}$ a constant depending only on the
parameter $\gamma$, and which may change from line to line.

\begin{lemma}
\label{LE1} Let $\phi\in\Phi$ and $\alpha\in(0,2]$. There is an increasing
function $\phi_{\alpha}:(0,1)\rightarrow R_{+}$ having $\lim_{x\downarrow
0}\phi_{\alpha}=0$ and $\phi\leq\phi_{\alpha}$ on $(0,1]$, such that for all
$G\in G(\phi,\alpha)$ and $\varepsilon\in(0,1]$ we have
\begin{equation}
\int_{\{|x|\leq\varepsilon\}}|x|^{q}F(dx)\leq\left\{
\begin{array}
[c]{ll}%
\frac{q}{q-\alpha}~\varepsilon^{q-\alpha}~\phi_{\alpha}(\varepsilon)\qquad &
\mbox{if }~q>\alpha\\
\phi_{\alpha}(\varepsilon)\qquad & \mbox{if }~q=\alpha=2,
\end{array}
\right.  \label{eq:Fmoments_near0}%
\end{equation}%
\begin{equation}
\int_{\{\varepsilon<|x|\leq1\}}|x|F(dx)\leq\left\{
\begin{array}
[c]{ll}%
\phi_{\alpha}(1) & \mbox{if }~\alpha<1\\
\phi_{\alpha}(\varepsilon)~\log(1/\varepsilon)\quad & \mbox{if }~\alpha=1\\
\phi_{\alpha}(\varepsilon)~\varepsilon^{1-\alpha} & \mbox{if }~\alpha>1.
\end{array}
\right.  \label{eq:Fmoments_awayfrom0}%
\end{equation}

\end{lemma}

\begin{proof}
First we define $\phi_{\alpha}$ as follows, for $x\in(0,1)$:
\[
\phi_{\alpha}(x)=\left\{
\begin{array}
[c]{ll}%
\frac{\phi(x)}{1-\alpha} & \mbox{if }~\alpha<1\\[2.5mm]%
\phi(x)+\frac{\phi(x)}{\sqrt{\log(1/x)}}+\phi\left(  1\wedge e^{-\sqrt
{\log(1/x)}}\right)  \qquad & \mbox{if }~\alpha=1\\[2.5mm]%
\phi(x)+\frac{\phi(\sqrt{x})}{\alpha-1}+\frac{\phi(1)}{\alpha-1}%
~x^{\frac{\alpha-1}{2}} & \mbox{if }~\alpha>1.
\end{array}
\right.
\]
It is clear that $\phi_{\alpha}(x)\rightarrow0$ as $x\downarrow0$, and that
$\phi\leq\phi_{\alpha}$ on $(0,1)$.

(\ref{eq:Fmoments_near0}) when $q=\alpha=2$ is trivial because $\phi\leq
\phi_{\alpha}$. When $q>\alpha$, Fubini Theorem and (\ref{eq:Gclass}) yield:
\begin{align*}
\int_{\{|x|\leq\varepsilon\}}|x|^{q}F(dx)  &  =\int_{\{|x|\leq\varepsilon
\}}F(dx)~q\int_{0}^{|x|}y^{q-1}dy~=~q\int_{0}^{\varepsilon}y^{q-1}%
F(|x|>y)~dy\\
&  \leq~q\int_{0}^{\varepsilon}\phi(y)y^{q-1-\alpha}dy~\leq~\frac{q}{q-\alpha
}~\phi(\varepsilon)~\varepsilon^{q-\alpha}%
\end{align*}
because $\phi$ is increasing: so we get (\ref{eq:Fmoments_near0}) again.

In a similar way, for every $z\in\lbrack\varepsilon,1]$ we get
\begin{align*}
\int_{\{\varepsilon<|x|\leq1\}}|x|F(dx)~~  &  =~\int_{\{\varepsilon
<|x|\leq1\}}F(dx)\int_{0}^{|x|}dy\\
&  =\int_{0}^{\varepsilon}F(\varepsilon<|x|\leq1)~dy+\int_{\varepsilon}%
^{z}F(y<|x|\leq1)~dy+\int_{z}^{1}F(y<|x|\leq1)~dy\\
&  \leq\phi(\varepsilon)\varepsilon^{1-\alpha}+\phi(z)\int_{\varepsilon}%
^{z}y^{-\alpha}~dy+\phi(1)\int_{z}^{1}y^{\alpha}~dy,
\end{align*}
Then in view of our definition of $\phi_{\alpha}$, a simple calculation allows
to deduce (\ref{eq:Fmoments_awayfrom0}), upon taking $z=1$ when $\alpha<1$,
and $z=1$ when $\alpha=1$ and $\varepsilon\geq1/e$, and $z=\exp-\sqrt
{\log(1/\varepsilon)}$ if $\alpha=1$ and $\varepsilon<1/e$, and $z=\sqrt
{\varepsilon}$ when $\alpha>1$.
\end{proof}

In view of (\ref{eq:Fmoments_near0}), for any pair $(G,\alpha)$ such that
$G\in G_{\alpha}$ we can introduce the following notation:
\begin{equation}
b^{\prime}(G,\alpha)=\left\{
\begin{array}
[c]{ll}%
b-\int_{\{|x|\leq1\}}xF(dx)\text{ \ \ \ } & \text{if \ }\alpha<1\\[2.5mm]%
b & \text{if \ }\alpha\geq1,
\end{array}
\qquad Z_{\Delta}(\alpha,\beta):=\Delta^{-1/\beta}\left(  Y_{\Delta}%
-b^{\prime}(G,\alpha)\Delta\right)  .\right.  \label{eq:bprime}%
\end{equation}
and we let $G_{\Delta,\alpha,\beta}^{\prime}$ denote the law of $Z_{\Delta
}(\alpha,\beta)$. Part (c) of the forthcoming lemma is not fully used here,
but will be in the companion paper on estimation.

\begin{lemma}
\label{LL2} a) If $G\in G_{\beta}$ then $G_{\Delta,\beta,\beta}^{\prime}$
converges to the Dirac mass $\delta_{0}$ as $\Delta\rightarrow0$.

b) If $\alpha\leq\beta$ and $\phi\in\Phi$ and $G^{n}$ is a sequence of
measures in $G(\phi,\alpha)$ and $\Delta_{n}\rightarrow0$, then the associated
sequence $G_{\Delta_{n},\alpha,\beta}^{\prime n}$ converges to the Dirac mass
$\delta_{0}$ as $n\rightarrow\infty$.

c) If $\alpha\leq\beta$ and $\phi\in\Phi$ there is a constant $C=C_{\alpha}$
such that for all functions $g$ with $|g(x)|\leq K(1\wedge|x|)$ and all
$\Delta\in(0,1]$, we have (with $\phi_{\alpha}$ like in the previous lemma):
\begin{equation}
G\in\mathcal{G}(\phi,\alpha)\qquad\Longrightarrow\qquad\mathbf{E}%
(|g(Z_{\Delta}(\alpha,\beta))|)\leq CK\Delta^{\frac{2(\beta-\alpha)}%
{\beta(2+\alpha)}}~\phi_{\alpha}\left(  \Delta^{\frac{2+\beta}{\beta
(2+\alpha)}}\right)  . \label{1}%
\end{equation}

\end{lemma}

\begin{proof}
Observe that (c)$\Rightarrow$(b)$\Rightarrow$(a), so we prove (c) only.

Let $\eta\in(0,1/2]$ to be chosen later. For any given $G\in G(K,\alpha)$ we
associate the L\'{e}vy process $Y$ and the characteristics $(b,0,F)$. Let
$F^{\prime}$ and $F^{\prime\prime}$ be the restrictions of $F$ to the sets
$[-\eta,\eta]$ and $[-\eta,\eta]^{c}$ respectively. We can decompose $Y$ into
the sum $Y_{t}=at+Y_{t}^{\prime}+Y_{t}^{\prime\prime}$, where $Y^{\prime}$ is
a L\'{e}vy process with characteristics $(0,0,F^{\prime})$ and $Y^{\prime
\prime}$ is a compound Poisson process with L\'{e}vy measure $F^{\prime\prime
}$, and $a=b-\int_{\{\eta<|x|\leq1\}}xF(dx)$. Then $a^{\prime}=a-b^{\prime
}(G,\alpha)$ is (recall (\ref{eq:bprime})):
\[
a^{\prime}=\left\{
\begin{array}
[c]{ll}%
\int_{\{|x|\leq\eta\}}xF(dx)\quad & \text{if \ }\alpha<1\\[2mm]%
-\int_{\{\eta<|x|\leq1\}}xF(dx)\quad & \text{if \ }\alpha\geq1.
\end{array}
\right.
\]
Therefore (\ref{eq:Fmoments_near0}) and (\ref{eq:Fmoments_awayfrom0}) yield
(for a constant $C=C_{\alpha}$ not depending on $G\in G(\phi,\alpha)$):
\begin{equation}
|a^{\prime}|\leq\left\{
\begin{array}
[c]{ll}%
C\eta^{1-\alpha}\phi_{\alpha}(\eta)\qquad & \text{if \ }\alpha\neq1\\[2.5mm]%
C\log(1/\eta)~\phi_{\alpha}(\eta)\qquad & \text{if \ }\alpha=1.
\end{array}
\right.  \label{60}%
\end{equation}

Also, since $Y^{\prime}$ has no drift, no Wiener part, and no jump bigger than
$1$, one knows (by differentiating (\ref{eq:Ychar}) for example) that
$E((Y_{t}^{\prime})^{2})=t\int x^{2}F^{\prime}(dx)$. Then
(\ref{eq:Fmoments_near0}) again yields for some $C=C_{\alpha}$:
\begin{equation}
\mathbf{E}(|Y_{\Delta}^{\prime}|^{2})\leq C\Delta~\eta^{2-\alpha}\phi_{\alpha
}(\eta). \label{61}%
\end{equation}

We set $Z_{\Delta}=Z_{\Delta}(\alpha,\beta)$. Since $|g|\leq K$ we have
$|g(Z_{\Delta})|\leq K$. If further $Y_{\Delta}^{\prime\prime}=0$, we have
also $Y_{\Delta}=a\Delta+Y_{\Delta}^{\prime}$, hence $Z_{\Delta}%
=\Delta^{-1/\beta}(Y_{\Delta}^{\prime}+a^{\prime}\Delta)$, hence
$|g(Z_{\Delta})|\leq K\Delta^{-1/\beta}(|Y_{\Delta}^{\prime}|+\Delta
|a^{\prime}|)$. Now, we have $P(Y_{\Delta}^{\prime\prime}\neq0)\leq\Delta
F^{\prime\prime}(I\!\!R)\leq\Delta\phi_{\alpha}(\eta)/\eta^{\alpha}$ because
$G\in G(\phi,\alpha)$. Therefore we deduce from (\ref{60}) and (\ref{61}) that
for some constant $C=C_{K,\alpha}$:
\[
\mathbf{E}(|g(Z_{\Delta})|)\leq\left\{
\begin{array}
[c]{ll}%
CK\left(  \Delta\eta^{-1}+\Delta^{1/2-1/\beta}\eta^{1/2}+\Delta^{1-1/\beta
}\log(1/\eta)\right)  \phi_{1}(\eta)\quad & \text{if \ }\alpha=1\\[2.5mm]%
CK\left(  \Delta\eta^{-\alpha}+\Delta^{1/2-1/\beta}\eta^{1-\alpha/2}%
+\Delta^{1-1/\beta}\eta^{1-\alpha}\right)  \phi_{\alpha}(\eta)\quad &
\text{otherwise}%
\end{array}
\right.
\]
as soon as $G\in G(\phi,\alpha)$. Then take $\eta=\Delta^{(2+\beta
)/\beta(2+\alpha)}$ to get (\ref{1}).
\end{proof}

\subsection{Fisher's information when $X=\sigma W+\theta Y$%
.\label{subsec:Fisher}}

>From independence of $W$ and $Y,$ the density of $X_{\Delta}$ in
(\ref{eq:X=sW+aY}) is the convolution (recall $G_{\Delta}=L(Y_{\Delta})$):
\begin{equation}
p_{\Delta}(x|\sigma,\beta,\theta,G)={\frac{1}{\sigma\Delta^{1/\beta}}}\int
G_{\Delta}(dy)h_{\beta}\left(  {\frac{x-\theta y}{\sigma\Delta^{1/\beta}}%
}\right)  . \label{eq:pdens_withY}%
\end{equation}

We now seek to characterize the entries of the full Fisher information matrix.
Since $h_{\beta}^{\prime}$ and $\breve{h}_{\beta}$ and $\dot{h}_{\beta}$ are
continuous and bounded, we can differentiate under the integral in
(\ref{eq:pdens_withY})\ to get
\begin{align}
\partial_{\sigma}p_{\Delta}(x|\sigma,\beta,\theta,G)  &  =-{\frac{1}%
{\sigma^{2}\Delta^{1/\beta}}}\int G_{\Delta}(dy)~\breve{h}_{\beta}\left(
{\frac{x-\theta y}{\sigma\Delta^{1/\beta}}}\right)  ,\label{eq:dpdsigma}\\
\partial_{\beta}p_{\Delta}(x|\sigma,\beta,\theta,G)  &  =v_{\Delta}%
(x|\sigma,\beta,\theta,G)-{\frac{\sigma\log\Delta}{\beta^{2}}}~\partial
_{\sigma}p_{\Delta}(x|\sigma,\beta,\theta,G),\label{eq:dpdbeta}\\
\partial_{\theta}p_{\Delta}(x|\sigma,\beta,\theta,G)  &  =-{\frac{1}%
{\sigma^{2}\Delta^{2/\beta}}}\int G_{\Delta}(dy)~y~h_{\beta}^{\prime}\left(
{\frac{x-\theta y}{\sigma\Delta^{1/\beta}}}\right)  , \label{eq:dpdtheta}%
\end{align}
where
\begin{equation}
v_{\Delta}(x|\sigma,\beta,\theta,G)={\frac{1}{\sigma\Delta^{1/\beta}}}\int
G_{\Delta}(dy)~\dot{h}_{\beta}\left(  {\frac{x-\theta y}{\sigma\Delta
^{1/\beta}}}\right)  . \label{eq:v}%
\end{equation}
The entries of the $(\sigma,\theta)$ block of the Fisher information matrix
are (leaving implicit the dependence on $(\sigma,\beta,\theta,G)$):
\begin{equation}
I_{\Delta}^{\sigma\sigma}=\int\frac{\partial_{\sigma}p_{\Delta}(x)^{2}%
}{p_{\Delta}(x)}~dx,\qquad I_{\Delta}^{\sigma\theta}=\int\frac{\partial
_{\sigma}p_{\Delta}(x)\partial_{\theta}p_{\Delta}(x)}{p_{\Delta}(x)}~dx,\qquad
I_{\Delta}^{\theta\theta}=\int\frac{\partial_{\theta}p_{\Delta}(x)^{2}%
}{p_{\Delta}(x)}~dx. \label{eq:I1}%
\end{equation}
When $\beta<2$, the other entries are
\begin{align}
I_{\Delta}^{\sigma\beta}  &  =J_{\Delta}^{\sigma\beta}-{\frac{\sigma\log
\Delta}{\beta^{2}}}~I_{\Delta}^{\sigma\sigma},\qquad I_{\Delta}^{\beta\theta
}=J_{\Delta}^{\beta\theta}-{\frac{\sigma\log\Delta}{\beta^{2}}}~I_{\Delta
}^{\sigma\theta},\label{eq:I2}\\
I_{\Delta}^{\beta\beta}  &  =J_{\Delta}^{\beta\beta}-{\frac{2\sigma\log\Delta
}{\beta^{2}}}~J_{\Delta}^{\sigma\beta}+{\frac{\sigma^{2}(\log\Delta)^{2}%
}{\beta^{4}}}~I_{\Delta}^{\sigma\sigma}, \label{eq:I3}%
\end{align}
where
\begin{equation}
J_{\Delta}^{\sigma\beta}=\int{\frac{\partial_{\sigma}p_{\Delta}(x)v_{\Delta
}(x)}{p_{\Delta}(x)}}~dx,\qquad J_{\Delta}^{\beta\beta}=\int{\frac{v_{\Delta
}(x)^{2}}{p_{\Delta}(x)}}~dx,\qquad J_{\Delta}^{\beta\theta}=\int
\frac{v_{\Delta}(x)\partial_{\theta}p_{\Delta}(x)}{p_{\Delta}(x)}~dx.
\label{eq:J}%
\end{equation}

\subsection{Proof of Theorem \ref{theorem:generalmajoration}.}

The proof is standard, and given for completeness, and given only in the case
where $\beta<2$ (when $\beta=2$ take $v=0$ below). What we need to prove is
that, for any $u,v\in R$, we have
\begin{align}
&  \int\frac{(u\partial_{\sigma}p_{\Delta}(x|\sigma,\beta,\theta
,G)+v\partial_{\beta}p_{\Delta}(x|\sigma,\beta,\theta,G))^{2}}{p_{\Delta
}(x|\sigma,\beta,\theta,G)}~dx\nonumber\\
&  \leq\int\frac{(u\partial_{\sigma}p_{\Delta}(x|\sigma,\beta,0,\delta
_{0})+v\partial_{\beta}p_{\Delta}(x|\sigma,\beta,0,\delta_{0}))^{2}}%
{p_{\Delta}(x|\sigma,\beta,0,\delta_{0})}~dx. \label{eq:majfi}%
\end{align}

We set
\[
q(x)=p_{\Delta}(x|\sigma,\beta,\theta,G),\qquad q_{0}(x)=p_{\Delta}%
(x|\sigma,\beta,0,\delta_{0}),
\]%
\[
r(x)=u\partial_{\sigma}p_{\Delta}(x|\sigma,\beta,\theta,G)+v\partial_{\beta
}p_{\Delta}(x|\sigma,\beta,\theta,G),\qquad r_{0}(x)=u\partial_{\sigma
}p_{\Delta}(x|\sigma,\beta,0,\delta_{0})+v\partial_{\beta}p_{\Delta}%
(x|\sigma,\beta,0,\delta_{0}).
\]
Observe that by (\ref{eq:pdens_withY}),
\[
q(x)=\int G_{\Delta}(dy)q_{0}(x-\theta y),
\]
hence
\[
r(x)=\int G_{\Delta}(dy)r_{0}(x-\theta y)
\]
as well. Apply Cauchy--Schwarz inequality to $G_{\Delta}$ with $r_{0}%
=\sqrt{q_{0}}~(r_{0}/\sqrt{q_{0}})~$ to get
\[
r(x)^{2}\leq q(x)\int G_{\Delta}(dy)~\frac{r_{0}(x-\theta y)^{2}}%
{q_{0}(x-\theta y)}.
\]
Then
\[
\int\frac{r(x)^{2}}{q(x)}~dx\leq\int dx\int G_{\Delta}(dy)~\frac
{r_{0}(x-\theta y)^{2}}{q_{0}(x-\theta y)}=\int\frac{r_{0}(z)^{2}}{q_{0}%
(z)}~dz
\]
by Fubini and a change of variable: this is exactly (\ref{eq:majfi}).

\subsection{Proof of Theorem \ref{theorem:sigma&beta}.}
Clearly (b) implies (a). If (b) fails, one can find $\phi\in\Phi$ and
$\alpha\in(0,\beta]$ and $\varepsilon>0$, and also a sequence $\Delta
_{n}\rightarrow0$ and a sequence $G^{n}$ of measures in $G(\phi,\alpha)$ and a
sequence of numbers $\theta_{n}$ converging to a limit $\theta$, such that
\[
\beta=2~\Rightarrow~\left\vert I_{\Delta_{n}}^{\sigma\sigma}(\sigma
,\beta,\theta_{n},G^{n})-{\frac{\mathcal{I}(\beta)}{\sigma^{2}}}\right\vert
\geq\varepsilon,
\]%
\[
\beta<2~\Rightarrow~\left\vert I_{\Delta_{n}}^{\sigma\sigma}(\sigma
,\beta,\theta_{n},G^{n})-{\frac{\mathcal{I}(\beta)}{\sigma^{2}}}\right\vert
+\left\vert \frac{I_{\Delta_{n}}^{\beta\beta}(\sigma,\beta,\theta_{n},G^{n}%
)}{(\log(1/\Delta))^{2}}-{\frac{\mathcal{I}(\beta)}{\beta^{4}}}\right\vert
+\left\vert \frac{I_{\Delta_{n}}^{\sigma\sigma}(\sigma,\beta,\theta_{n}%
,G^{n})}{\log(1/\Delta)}-{\frac{\mathcal{I}(\beta)}{\sigma\beta^{2}}%
}\right\vert \geq\varepsilon
\]
for all $n$.

In other words, to prove (a) and (b) it is enough to prove the following: let
$\phi\in\Phi$ and $\alpha\in(0,\beta]$ and $\Delta_{n}\rightarrow0$ and
$\theta_{n}\rightarrow\theta$ and $G^{n}$ be a sequence in $G(\phi,\alpha)$;
then we have:
\begin{equation}
I_{\Delta_{n}}^{\sigma\sigma}(\sigma,\beta,\theta_{n},G^{n})\rightarrow
{\frac{\mathcal{I}(\beta)}{\sigma^{2}}}, \label{eq:conver}%
\end{equation}%
\begin{equation}
\beta<2\quad\Rightarrow\quad\frac{I_{\Delta_{n}}^{\beta\beta}(\sigma
,\beta,\theta_{n},G^{n})}{(\log(1/\Delta))^{2}}\rightarrow{\frac
{\mathcal{I}(\beta)}{\beta^{4}}},\qquad\frac{I_{\Delta_{n}}^{\sigma\beta
}(\sigma,\beta,\theta_{n},G^{n})}{\log(1/\Delta)}\rightarrow{\frac
{\mathcal{I}(\beta)}{\sigma\beta^{2}}}. \label{eq:conver2}%
\end{equation}
\medskip

Let us proceed to proving (\ref{eq:conver}). The change of variable
$x\leftrightarrow(x-\theta b^{\prime}(G^{n},\alpha))/\sigma\Delta^{1/\beta}$
in (\ref{eq:I1}) leads to $I_{\Delta}^{\sigma\sigma}(\sigma,\beta
,\theta,G)=\frac{1}{\sigma^{2}}\int s_{\Delta,\theta,G}(x)dx$, where
\begin{equation}
s_{\Delta,\theta,G}(x)=\frac{\left(  \int G_{\Delta,\alpha,\beta}^{\prime
}(du)~\breve{h}_{\beta}(x-u\theta/\sigma)\right)  ^{2}}{\int G_{\Delta
,\alpha,\beta}^{\prime}(du)h_{\beta}(x-u\theta/\sigma)}. \label{eq:s}%
\end{equation}
Since $h_{\beta}$ and $\breve{h}_{\beta}$ are continuous and bounded, we
deduce from Lemma \ref{LL2} that
\[
\int G^{\prime n}n_{\Delta_{n},\alpha,\beta}(du)\breve{h}_{\beta}%
(x-u\theta_{n}/\sigma)\rightarrow\breve{h}_{\beta}(x),\qquad\int G^{\prime
n}n_{\Delta_{n},\alpha,\beta}(du)h(x-u\theta_{n}/\sigma)\rightarrow h(x).
\]
Thus $s_{\Delta_{n},\theta_{n},G^{n}}(x)\rightarrow\widetilde{h}_{\beta}(x)$
for all $x$, and Fatou's Lemma yields
\[
\liminf_{n}I_{\Delta_{n}}^{\sigma\sigma}(\sigma,\beta,\theta_{n},G^{n}%
)\geq\frac{1}{\sigma^{2}}\int\widetilde{h}_{\beta}(x)dx=\frac{\mathcal{I}%
(\beta)}{\sigma^{2}}.
\]
This, combined with (\ref{eq:Icallig}) and (\ref{eq:Iss_Y=0}) and
(\ref{FI-MAJ}), gives (\ref{eq:conver}). \medskip

Now suppose that $\beta<2$ and recall (\ref{eq:I2}) and (\ref{eq:I3}): with
the notation $J_{\Delta}^{\sigma\beta}(\sigma,\beta,\theta,G)$, etc... of
(\ref{eq:J}), we see that
\begin{equation}
J_{\Delta}^{\sigma\beta}(\sigma,\beta,\theta,G)\leq\sqrt{I_{\Delta}%
^{\sigma\sigma}(\sigma,\beta,\theta,G)J_{\Delta}^{\beta\beta}(\sigma
,\beta,\theta,G)} \label{eq:CS2}%
\end{equation}
by a first application of Cauchy--Schwarz inequality. A second application of
the same plus (\ref{eq:pdens_withY}) and (\ref{eq:v}) yield
\[
v_{\Delta}(x|\sigma,\beta,\theta,G)^{2}\leq{\frac{p_{\Delta}(x|\sigma
,\beta,\theta,G)}{\sigma\Delta^{1/\beta}}}\int G_{\Delta}(dy)~\frac{\dot
{h}_{\beta}^{2}}{h_{\beta}}\left(  {\frac{x-\theta y}{\sigma\Delta^{1/\beta}}%
}\right)  .
\]
Then Fubini and the change of variable $x\leftrightarrow(x-\theta
y)/\sigma\Delta^{1/\beta}$ in (\ref{eq:J}) leads to
\begin{equation}
J_{\Delta}^{\beta\beta}(\sigma,\beta,\theta,G)\leq\mathcal{K}(\beta).
\label{eq:majJbeta}%
\end{equation}
Then (\ref{eq:conver2}) readily follow from (\ref{eq:conver}), (\ref{eq:CS2}),
(\ref{eq:majJbeta}) and also (\ref{eq:I2}) and (\ref{eq:I3}). \bigskip

It remains to prove (c). If we put together the majorations (\ref{FI-MAJ}),
(\ref{eq:CS2}) and (\ref{eq:majJbeta}) and also (\ref{eq:I2}) and
(\ref{eq:I3}), we see that it is enough to prove that
\begin{equation}
\limsup_{n}I_{\Delta_{n}}^{\sigma\sigma}(\sigma,\beta,\theta,G^{n}%
)<\frac{\mathcal{I}(\beta)}{\sigma^{2}}. \label{eq:limsup}%
\end{equation}

Let $\rho_{n}=\Delta_{n}^{1/\alpha_{n}-1/\beta}$, which by our assumption on
$\Delta_{n}$ goes to $1$. The measure $G_{\Delta_{n},\alpha_{n},\beta}^{\prime
n}$ admits the density $x\mapsto g_{n}(x)=h_{\alpha_{n}}(x\rho_{n})/\rho_{n}$,
which converges to $h_{\beta}(x)$; so $G_{\Delta_{n},\alpha_{n},\beta}^{\prime
n}$ weakly converges to the stable law with density $h_{\beta}$. Then, exactly
as in the previous proof, we get that
\begin{equation}
s_{\Delta_{n},\theta,G^{n}}(x)\rightarrow s(x):={\frac{\left(  \int h_{\beta
}(u)\breve{h}_{\beta}(x-u\theta/\sigma)du\right)  ^{2}}{\int h_{\beta
}(u)h_{\beta}(x-u\theta/\sigma)~du}}. \label{100}%
\end{equation}
On the other hand $|\widetilde{h}_{\beta}(y)|~\leq C(1\wedge1/|y|^{1+\beta})$
for some constant $C$, and also $g_{n}^{\prime}(y)\leq C(1\wedge
1/|y|^{1+\alpha_{n}})$ with $C$ not depending on $n$. Using once more
Cauchy--Schwarz, we deduce from (\ref{eq:s}) that
\begin{align}
s_{\sigma,\Delta_{n},G^{n}}(x)  &  \leq\int g_{n}(u)\widetilde{h}_{\beta
}(x-u\theta/\sigma)~du\nonumber\\
&  \leq s^{\prime}(x):=C\int\left(  1\wedge{\frac{1}{|u|^{1+\beta-\varepsilon
}}}\right)  \left(  1\wedge{\frac{1}{|x-u\theta/\sigma|^{1+\beta}}}\right)
~du. \label{eq:CS}%
\end{align}
for still another constant $C$, as soon as $\alpha_{n}>\beta-\varepsilon$ for
some fixed $\varepsilon\in(0,\beta)$. But $\int s^{\prime}(x)dx$ is finite, so
(\ref{100}) and the dominated convergence theorem yield
\begin{equation}
I_{\Delta_{n}}^{\sigma\sigma}(\sigma,\beta,\theta,G^{n})\rightarrow\frac
{1}{\sigma^{2}}\int s(x)dx. \label{FI10}%
\end{equation}

Finally, exactly as for (\ref{eq:CS}) we deduce from the Cauchy--Schwarz
inequality and from the fact that the functions $\sqrt{h}_{\beta}$ and
$\breve{h}_{\beta}/\sqrt{h}_{\beta}$ are not Lebesgue--almost surely multiple
one from the other, while $h_{\beta}>0$ identically, that in fact $s(x)<\int
h_{\beta}(u)\widetilde{h}_{\beta}(x-u\theta/\sigma)du$ for all $x$. Therefore
\[
\int s(x)dx<\int dx\int h_{\beta}(u)\widetilde{h}_{\beta}(x-u\theta
/\sigma)du=\int h_{\beta}(u)du\int\widetilde{h}_{\beta}(y)dy=\int\widetilde
{h}_{\beta}(y)dy=\mathcal{I}(\beta),
\]
and (\ref{FI10}) yields that (\ref{eq:limsup}) holds, hence (c).

\subsection{Proof of Theorem \ref{theorem:theta}.}

Cauchy--Schwarz inequality gives us, by (\ref{eq:pdens_withY}) and
(\ref{eq:dpdtheta}):
\[
|\partial_{\theta}p_{\Delta}(x|\Delta,\beta,\theta,G)|^{2}\leq\frac{1}%
{\sigma^{3}\Delta^{3/\beta}}~p_{\Delta}(x|\Delta,\beta,\theta,G)~\int
G_{\Delta}(dy)~y^{2}~\widetilde{h}_{\beta}\left(  \frac{x-\theta y}%
{\sigma\Delta^{1:\beta}}\right)  .
\]
Plugging this into (\ref{eq:I1}), applying Fubini and doing the change of
variable $x\leftrightarrow(x-\theta y)/\sigma\Delta^{1/\beta}$ leads to
\[
I_{\Delta}^{\theta\theta}\leq\frac{1}{\sigma^{2}\Delta^{2/\beta}}~\int
G_{\Delta}(dy)y^{2}~\int\frac{h_{\beta}^{\prime}(x)^{2}}{h_{\beta}(x)}~dx.
\]
Since $E(Y_{\Delta}^{2})=m\Delta^{2}+\delta\Delta$, we readily deduce
(\ref{FI-INtheta1}).

\subsection{Proof of Theorem \ref{theorem:S+D}.}

When $Y_{t}=t$ we gave $G_{\Delta}=\delta_{\Delta}$. Then (\ref{eq:Fisher_S+D}%
) follows directly from applying the formulae (\ref{eq:dpdsigma}) and
(\ref{eq:dpdtheta}) and from the change of variable $x\leftrightarrow
(x-\theta)/\sigma\Delta^{1/\beta}$ in (\ref{eq:I1}), after observing that the
function $\breve{h}_{\beta}h_{\beta}^{\prime}/h_{\beta}$ is integrable and
odd, hence has a vanishing Lebesgue integral.

\subsection{Proof of Theorem \ref{theorem:S+P}.}

a) We first introduce some notation to be used also for the proof of Theorem
\ref{theorem:S+CP}. We suppose that $Y$ is a compound Poisson process with
arrival rate $\lambda$ and law of jumps $\mu$, and $\mu_{k}$ is the $k$th fold
convolution of $\mu$. So we have
\[
G_{\Delta}=e^{-\lambda\Delta}\sum_{k=0}^{\infty}\frac{(\lambda\Delta)^{k}}%
{k!}~\mu_{k}.
\]
Set
\begin{align}
\gamma_{\Delta}^{(1)}(k,x)  &  =\frac{1}{\sigma\Delta^{1/\beta}}\int\mu
_{k}(du)~h_{\beta}\left(  \frac{x-\theta u}{\sigma\Delta^{1/\beta}}\right)
,\label{eq:gamma01}\\
\gamma_{\Delta}^{(2)}(k,x)  &  =\frac{1}{\sigma^{2}\Delta^{1/\beta}}\int
\mu_{k}(du)~\breve{h}_{\beta}\left(  \frac{x-\theta u}{\sigma\Delta^{1/\beta}%
}\right)  ,\label{eq:gamma02}\\
\gamma_{\Delta}^{(3)}(k,x)  &  =\frac{1}{\sigma^{2}\Delta^{2/\beta}}\int
\mu_{k}(du)~u~h_{\beta}^{\prime}\left(  \frac{x-\theta u}{\sigma
\Delta^{1/\beta}}\right)  . \label{eq:gamma03}%
\end{align}
We have (recall that $\mu_{0}=\delta_{0}$):
\begin{align}
p_{\Delta}(x|\sigma,\beta,\theta,G)  &  =e^{-\lambda\Delta}\sum_{k=0}^{\infty
}\frac{(\lambda\Delta)^{k}}{k!}~\gamma_{\Delta}^{(1)}(k,x),\label{eq:denCP}\\
\partial_{\sigma}p_{\Delta}(x|\sigma,\beta,\theta,G)  &  =-e^{-\lambda\Delta
}\sum_{k=0}^{\infty}\frac{(\lambda\Delta)^{k}}{k!}~\gamma_{\Delta}%
^{(2)}(k,x),\label{eq:sidenCP}\\
\partial_{\theta}p_{\Delta}(x|\sigma,\beta,\theta,G)  &  =-e^{-\lambda\Delta
}\sum_{k=1}^{\infty}\frac{(\lambda\Delta)^{k}}{k!}~\gamma_{\Delta}^{(3)}(k,x),
\label{eq:tedenCP}%
\end{align}

Omitting the mention of $(\sigma,\beta,\theta,G)$, we also set
\begin{equation}
i=2,3:\quad\Gamma_{\Delta}^{(i)}(k,k^{\prime})=\int\frac{\gamma_{\Delta}%
^{(i)}(k,x)\gamma_{\Delta}^{(i)}(k^{\prime},x)}{p_{\Delta}(x)}~dx,\qquad
\Gamma_{\Delta}^{(4)}(k,k^{\prime})=\int\frac{\gamma_{\Delta}^{(2)}%
(k,x)\gamma_{\Delta}^{(3)}(k^{\prime},x)}{p_{\Delta}(x)}~dx. \label{eq:Gamma1}%
\end{equation}
By Cauchy--Schwarz inequality, we have
\begin{equation}
i=2,3:\quad\left\vert \Gamma_{\Delta}^{(i)}(k,k^{\prime})\right\vert \leq
\sqrt{\Gamma_{\Delta}^{(i)}(k,k)\Gamma_{\Delta}^{(i)}(k^{\prime},k^{\prime}%
)},\qquad\left\vert \Gamma_{\Delta}^{(4)}(k,k^{\prime})\right\vert \leq
\sqrt{\Gamma_{\Delta}^{(2)}(k,k)\Gamma_{\Delta}^{(3)}(k^{\prime},k^{\prime})}.
\label{eq:CSGamma}%
\end{equation}
For any $k\geq0$ we have $p_{\Delta}(x)\geq e^{-\lambda\Delta}\frac
{(\lambda\Delta)^{k}}{k!}~\gamma_{\Delta}^{(1)}(k,x)$. Therefore
\begin{equation}
i=2,3:\qquad\Gamma_{\Delta}^{(i)}(k,k)\leq\frac{e^{\lambda\Delta}~k!}%
{(\lambda\Delta)^{k}}~\int\frac{\gamma_{\Delta}^{(i)}(k,x)^{2}}{\gamma
_{\Delta}^{(1)}(k,x)}~dx. \label{eq:MajorationGamma}%
\end{equation}

Finally, if we plug (\ref{eq:sidenCP}) and (\ref{eq:tedenCP}) into
(\ref{eq:I1}), we get
\begin{align}
I_{\Delta}^{\sigma\theta}  &  =e^{-2\lambda\Delta}\sum_{k=0}^{\infty}%
\sum_{l=1}^{\infty}\frac{(\lambda\Delta)^{k+l}}{k!~l!}~\Gamma_{\Delta}%
^{(4)}(k,l),\label{eq:sigthep}\\
I_{\Delta}^{\theta\theta}  &  =e^{-2\lambda\Delta}\sum_{k=1}^{\infty}%
\sum_{l=1}^{\infty}\frac{(\lambda\Delta)^{k+l}}{k!~l!}~\Gamma_{\Delta}%
^{(3)}(k,l). \label{eq:thethep}%
\end{align}
\medskip

b) Now we can proceed to the proof of Theorem \ref{theorem:S+P}. When $Y$ is a
standard Poisson process, we have $\lambda=1$ and $\mu_{k}=\varepsilon_{k}$.
Therefore we get
\begin{align}
\gamma_{\Delta}^{(1)}(k,x)  &  =\frac{1}{\sigma\Delta^{1/\beta}}~h_{\beta
}\left(  \frac{x-\theta k}{\sigma\Delta^{1/\beta}}\right)  , \label{eq:gamma1}%
\\
\gamma_{\Delta}^{(2)}(k,x)  &  =\frac{1}{\sigma^{2}\Delta^{1/\beta}}~\breve
{h}_{\beta}\left(  \frac{x-\theta k}{\sigma\Delta^{1/\beta}}\right)
,\label{eq:gamma2}\\
\gamma_{\Delta}^{(3)}(k,x)  &  =\frac{k}{\sigma^{2}\Delta^{2/\beta}}~h_{\beta
}^{\prime}\left(  \frac{x-\theta k}{\sigma\Delta^{1/\beta}}\right)  .
\label{eq:gamma3}%
\end{align}
Plugging this into (\ref{eq:MajorationGamma}) yields
\begin{equation}
\Gamma_{\Delta}^{(2)}(k,k)\leq\frac{e^{\Delta}~k!}{\sigma^{2}\Delta^{k}%
}~\mathcal{I}(\beta),\qquad\Gamma_{\Delta}^{(3)}(k,k)\leq\frac{e^{\Delta
}~k^{2}~k!}{\sigma^{2}\Delta^{k+2/\beta}}~\mathcal{J}(\beta).
\label{eq:majorationGammaP}%
\end{equation}

Recall that (\ref{eq:Iss_Y=Poisson}) follows from Theorem
\ref{theorem:sigma&beta}, so we need to prove (\ref{eq:Isa_Y=Poisson}) and
(\ref{eq:Iaa_Y=Poisson}). In view of (\ref{eq:sigthep}) and (\ref{eq:thethep}%
), this amounts to proving the following two properties:
\[
\sum_{k=0}^{\infty}\sum_{l=1}^{\infty}\frac{\Delta^{k+l+1/\beta-1/2}}%
{k!~l!}~\Gamma_{\Delta}^{(4)}(k,l)\rightarrow0,\qquad\sum_{k=1}^{\infty}%
\sum_{l=1}^{\infty}\frac{\Delta^{k+l+2/\beta-1}}{k!~l!}~\Gamma_{\Delta}%
^{(3)}(k,l)\rightarrow\frac{1}{\sigma^{2}}~\mathcal{J}(\beta).
\]
If we use (\ref{eq:CSGamma}) and (\ref{eq:majorationGammaP}), it is easily
seen that the sum of all summands in the first (resp. second) left side above,
except the one for $k=0$ and $l=1$ (resp. $k=l=1$) goes to $0$. So we are left
to prove
\begin{equation}
\Delta^{1/2+1/\beta}~\Gamma_{\Delta}^{(4)}(0,1)\rightarrow0,\qquad
\Delta^{1+2/\beta}~\Gamma_{\Delta}^{(3)}(1,1)\rightarrow\frac{1}{\sigma^{2}%
}~\mathcal{J}(\beta). \label{eq:limits}%
\end{equation}

Let $\omega=\frac{1}{2(1+\beta)}$, so $(1+\frac{1}{\beta})(1-\omega
)=1+1/2\beta$. Assume first that $\beta<2$. Then if $|x|\leq|\theta
/\sigma\Delta^{1/\beta}|^{\omega}$, we have for some constant $C\in(0,\infty
)$, possibly depending on $\theta$, $\sigma$ and $\beta$, and which changes
from an occurrence to the other, and provided $\Delta\leq(2\theta
/\sigma)^{\beta}$:
\[
h_{\beta}(x)\geq C\Delta^{\omega(1+1/\beta)},\qquad h_{\beta}(x+r\theta
/\sigma\Delta^{1/\beta})\leq C\Delta^{1+1/\beta}%
\]
when $r\in Z\backslash\{0\}$, and thus
\[
|x|\leq|\theta/\sigma\Delta^{1/\beta}|^{\omega}\quad\Rightarrow\quad h_{\beta
}(x+r\theta/\sigma\Delta^{1/\beta})\leq Ch_{\beta}(x)\Delta^{(1+1/\beta
)(1-\omega)}=Ch_{\beta}(x)\Delta^{1+1/2\beta}.
\]
When $\beta=2$, a simple computation on the normal density shows that the
above property also holds. Therefore in view of (\ref{eq:denCP}) and
(\ref{eq:gamma1}) we deduce that in all cases,
\begin{align}
\left\vert \frac{x-\theta}{\sigma\Delta^{1/\beta}}\right\vert \leq\left(
\frac{\theta}{\sigma\Delta^{1/\beta}}\right)  ^{\omega}\quad\Rightarrow\quad
p_{\Delta}(x)  &  \leq\frac{e^{-\Delta}}{\sigma\Delta^{1/\beta}}~h_{\beta
}\left(  \frac{x-\theta}{\sigma\Delta^{1/\beta}}\right)  \left(
\Delta+C\Delta^{1+\beta/2}\left(  1+\sum_{k=2}^{\infty}\frac{\Delta^{k}}%
{k!}\right)  \right) \nonumber\\
&  \leq\frac{e^{-\Delta}}{\sigma\Delta^{1/\beta-1}}~h_{\beta}\left(
\frac{x-\theta}{\sigma\Delta^{1/\beta}}\right)  \left(  1+C\Delta^{\beta
/2}\right)  .\nonumber
\end{align}
By (\ref{eq:Gamma1}) it follows that
\begin{align}
\Gamma_{\Delta}^{(3)}(1,1)  &  \geq\frac{e^{\Delta}}{\sigma^{3}\Delta
^{1+3/\beta}}\int_{\{|(x-\theta)/\sigma\Delta^{1/\beta}|\leq(\theta
/\sigma\Delta^{1/\beta})^{\omega}\}}\frac{h_{\beta}^{\prime}((x-\theta
)/\sigma\Delta^{1/\beta})^{2}}{h_{\beta}((x-\theta)/\sigma\Delta^{1/\beta}%
)}~dx\nonumber\\
&  \geq\frac{e^{\Delta}}{\sigma^{2}\Delta^{1+2/\beta}}\int_{\{|x|\leq
(\theta/\sigma\Delta^{1/\beta})^{\omega}\}}\frac{h_{\beta}^{\prime}(x)^{2}%
}{h_{\beta}(x)}~dx\nonumber
\end{align}
We readily deduce that $\liminf_{\Delta\rightarrow0}\Delta^{1+2/\beta}%
\Gamma_{\Delta}^{(3)}(1,1)\geq J(\beta)/\sigma^{2}$. On the other hand,
(\ref{eq:majorationGammaP}) yields $\limsup_{\Delta\rightarrow0}%
\Delta^{1+2/\beta}\Gamma_{\Delta}^{(3)}(1,1)\leq J(\beta)/\sigma^{2}$, and
thus the second part of (\ref{eq:limits}) is proved.

Finally $\breve{h}_{\beta}/h_{\beta}$ is bounded, so (\ref{eq:gamma2}) and
(\ref{eq:gamma3}) and the fact that $p_{\Delta}(x)\geq e^{-\Delta}h_{\beta
}(x/\sigma\Delta^{1/\beta})/\sigma\Delta^{1/\beta}$ yield
\[
\left\vert \Gamma_{\Delta}^{(4)}(0,1)\right\vert \leq\frac{e^{\Delta}}%
{\sigma^{3}\Delta^{2/\beta}}\int\left\vert h_{\beta}^{\prime}\left(
\frac{x-\theta}{\sigma\Delta^{1/\beta}}\right)  \right\vert ~dx\leq
\frac{e^{\Delta}}{\sigma^{2}\Delta^{1/\beta}}\int|h_{\beta}^{\prime}(x)|~dx,
\]
and the first part of (\ref{eq:limits}) readily follows.

\subsection{Proof of Theorem \ref{theorem:S+CP}.}

We use the same notation than in the previous proof, but here the measure
$\mu_{k}$ has a density $f_{k}$ for all $k\geq1$, which further is
differentiable and satisfies (\ref{eq:cond-f}) uniformly in $k$, while we
still have $\mu_{0}=\delta_{0}$. Exactly as in (\ref{eq:hfunctions}), we set
\[
\breve{f}_{k}(u)=uf_{k}^{\prime}(u)+f_{k}(u).
\]
Recall the Fisher information $\mathcal{L}$ defined in (\ref{eq:FishMult}),
which corresponds to estimating $\theta$ in a model where one observes a
variable $\theta U$, with $U$ having the law $\mu$. Now if we have $n$
independent variables $U_{i}$ with the same law $\mu$, the Fisher information
associated with the observation of $\theta U_{i}$ for $i=1,\ldots,n$ is of
course $n\mathcal{L}$, and if instead one observes only $\theta(U_{1}%
+\ldots+U_{n})$, one gets a smaller Fisher information $\mathcal{L}_{n}\leq
n\mathcal{L}$. In other words, we have
\begin{equation}
\mathcal{L}_{n}:=\int\frac{\breve{f}_{n}(u)^{2}}{f_{n}(u)}~du\leq
n\mathcal{L}. \label{eq:FishMultn}%
\end{equation}

Taking advantage of the fact that $\mu_{k}$ has a density, for all $k\geq1$ we
can rewrite $\gamma_{\Delta}^{(i)}(k,x)$ as follows (using further an
integration by parts when $i=3$ and the fact that each $f_{k}$ satisfies
(\ref{eq:cond-f})):%

\begin{align}
\gamma_{\Delta}^{(1)}(k,x)  &  =\frac{1}{\theta}\int h_{\beta}(y)f_{k}\left(
\frac{x-y\sigma\Delta^{1/\beta}}{\theta}\right)  ~dy\label{eq:gammaPC1}\\
\gamma_{\Delta}^{(2)}(k,x)  &  =\frac{1}{\sigma\theta}\int\breve{h}_{\beta
}(y)f_{k}\left(  \frac{x-y\sigma\Delta^{1/\beta}}{\theta}\right)
~dy\label{eq:gammaPC2}\\
\gamma_{\Delta}^{(3)}(k,x)  &  =\frac{1}{\theta^{2}}\int h_{\beta}(y)\breve
{f}_{k}\left(  \frac{x-y\sigma\Delta^{1/\beta}}{\theta}\right)  ~dy.
\label{eq:gammaPC3}%
\end{align}
Since the $f_{k}$'s satisfy (\ref{eq:cond-f}) uniformly in $k$, we readily
deduce that
\begin{equation}
k\geq1,~i=1,2,3\quad\Rightarrow\quad|\gamma_{\Delta}^{(i)}(k,x)|\leq
C,\qquad\gamma_{\Delta}^{(1)}(1,x)\rightarrow\frac{1}{\theta}~f\left(
\frac{x}{\theta}\right)  ,\quad\gamma_{\Delta}^{(3)}(1,x)\rightarrow\frac
{1}{\theta^{2}}~\breve{f}\left(  \frac{x}{\theta}\right)  ,
\label{eq:majgamma}%
\end{equation}

Let us start with the lower bound. Since $\gamma_{\Delta}^{(1)}(0,x)=\frac
{1}{\sigma\Delta^{1/\beta}}~h_{\beta}\left(  \frac{x}{\sigma\Delta^{1/\beta}%
}\right)  $, we deduce from (\ref{eq:dh^n}) and (\ref{eq:majgamma}) and
(\ref{eq:denCP}) that
\begin{equation}
\frac{1}{\Delta}~p_{\Delta}(x)\rightarrow\frac{c_{\beta}~\sigma^{\beta}%
}{|x|^{1+\beta}}+\frac{\lambda}{\theta}~f\left(  \frac{x}{\theta}\right)
\label{AAA}%
\end{equation}
as soon as $x\neq0$, and with the convention $c_{2}=0$. In a similar way, we
deduce from (\ref{eq:majgamma}) and (\ref{eq:tedenCP}):
\begin{equation}
\frac{1}{\Delta}~\partial_{\theta}p_{\Delta}(x)\rightarrow-\frac{\lambda
}{\theta^{2}}~\breve{f}\left(  \frac{x}{\theta}\right)  . \label{BBB}%
\end{equation}
Then plugging (\ref{AAA}) and (\ref{BBB}) into the last equation in
(\ref{eq:I1}), we conclude by Fatou's Lemma and after a change of variable,
that
\[
\liminf_{\Delta\rightarrow0}~\frac{I_{\Delta}^{\theta\theta}}{\Delta}\geq
\frac{\lambda^{2}}{\theta^{2}}~\int\frac{\breve{f}(x)^{2}}{\lambda
f(x)+c_{\beta}~\sigma^{\beta}/\theta^{\beta}|x|^{1+\beta}}~dx.
\]
\medskip

It remains to prove (\ref{eq:Isa_Y=CPoisson}) and the upper bound in
(\ref{eq:Iaa_Y=CPoisson}) (including when $\beta=2$). By Cauchy--Schwarz, we
get (using successively the two equivalent versions for $\gamma_{\Delta}%
^{(i)}(k,x)$):
\[
\gamma_{\Delta}^{(2)}(k,x)^{2}\leq\frac{1}{\sigma^{3}}~\gamma_{\Delta}%
^{(1)}(k,x)\int\mu_{k}(du)~\widetilde{h}_{\beta}((x-\theta u)/\sigma
\Delta^{1/\beta}).
\]%
\[
\gamma_{\Delta}^{(3)}(k,x)^{2}\leq\frac{1}{\theta^{3}}~\gamma_{\Delta}%
^{(1)}(k,x)\int h_{\beta}(y)~\frac{\breve{f}_{k}((x-y\sigma\Delta^{1/\beta
})/\theta)^{2}}{f_{k}((x-y\sigma\Delta^{1/\beta})/\theta)}~dy.
\]
Then it follows from (\ref{eq:MajorationGamma}) and (\ref{eq:FishMultn}) that
\begin{equation}
\Gamma_{\Delta}^{(2)}(k,k)\leq\frac{e^{\lambda\Delta}~k!}{\sigma^{2}%
(\lambda\Delta)^{k}}~\mathcal{I}(\beta),\qquad\Gamma_{\Delta}^{(3)}%
(k,k)\leq\frac{e^{\lambda\Delta}~k~k!}{\theta^{2}(\lambda\Delta)^{k}%
}~\mathcal{L}. \label{eq:MajGamma}%
\end{equation}

We also need an estimate for $\Gamma_{\Delta}^{(4)}(0,1)$. By
(\ref{eq:gamma01}) and (\ref{eq:gamma02}) and (\ref{eq:denCP}) we obtain
$p_{\Delta}(x)\geq e^{-\lambda\Delta}~h_{\beta}(x/\sigma\Delta^{1/\beta
})/\sigma\Delta^{1/\beta}$ and $\gamma_{\Delta}^{(2)}(0,x)=\breve{h}_{\beta
}(x/\sigma\Delta^{1/\beta})/\sigma^{2}\Delta^{1/\beta}$. Then use
(\ref{eq:gammaPC3}) and the definition (\ref{eq:Gamma1}) to get
\begin{align}
\left\vert \Gamma_{\Delta}^{(4)}(0,1)\right\vert  &  \leq\frac{1}{\sigma
\theta^{2}}\int\left\vert \frac{\breve{h}_{\beta}}{h_{\beta}}\left(  \frac
{x}{\sigma\Delta^{1/\beta}}\right)  \right\vert ~dx\int h_{\beta}(y)\left\vert
\breve{f}\left(  \frac{x-y\sigma\Delta^{1/\beta}}{\theta}\right)  \right\vert
~dy\nonumber\\
&  =\frac{1}{\sigma\theta^{2}}\int h_{\beta}(y)~dy\int\left\vert \frac
{\breve{h}_{\beta}}{h_{\beta}}\left(  \frac{x}{\sigma\Delta^{1/\beta}}\right)
~\breve{f}\left(  \frac{x-y\sigma\Delta^{1/\beta}}{\theta}\right)  \right\vert
~dx~\leq~C, \label{eq:majGamma4}%
\end{align}
where the last inequality comes from the facts that $\breve{h}_{\beta
}/h_{\beta}$ is bounded and that $\breve{f}$ is integrable (due to
(\ref{eq:cond-f})).

At this stage we use (\ref{eq:CSGamma}), together with (\ref{eq:sigthep}) and
(\ref{eq:thethep}) and the fact that $2|xy|\leq ax^{2}+y^{2}/a$ for all $a>0$.
Taking arbitrary constants $a_{kl}>0$, we deduce from (\ref{eq:MajGamma})
that
\begin{align}
I_{\Delta}^{\theta\theta}  &  \leq\frac{\mathcal{L}}{\theta^{2}}%
~e^{-\lambda\Delta}\left(  \sum_{k=1}^{\infty}\frac{(\lambda\Delta)^{k}~k}%
{k!}+2\sum_{k=1}^{\infty}\sum_{l=k+1}^{\infty}\frac{(\lambda\Delta)^{k+l}%
}{k!~l!}~\left(  \frac{k~k!}{(\lambda\Delta)^{k}}\frac{l~l!}{(\lambda
\Delta)^{l}}\right)  ^{1/2}\right) \nonumber\\
&  \leq\frac{\mathcal{L}}{\theta^{2}}~\left(  \lambda\Delta+\sum_{k=1}%
^{\infty}\sum_{l=k+1}^{\infty}\frac{(\lambda\Delta)^{l}~k}{l!}~a_{kl}%
+\sum_{k=1}^{\infty}\sum_{l=k+1}^{\infty}\frac{(\lambda\Delta)^{k}~l}%
{k!}~\frac{1}{a_{kl}}\right) \nonumber
\end{align}
Then if we take $a_{kl}=(\lambda\Delta)^{\frac{k-l}{2}}$ for $l>k$, a simple
computation shows that indeed
\[
I_{\Delta}^{\theta\theta}\leq\frac{\mathcal{L}}{\theta^{2}}~\left(
\lambda\Delta+C\Delta^{3/2}\right)
\]
for some constant $C$, and thus we get the upper bound in
(\ref{eq:Iaa_Y=CPoisson}). In a similar way, and replacing $L/\theta^{2}$
above by the supremum between $L/\theta^{2}$ and $I(\beta)/\sigma^{2}$, we see
that in (\ref{eq:sigthep}) the sum of the absolute values of all summands
except the one for $k=0$ and $l=1$ is smaller than a constant times
$\Lambda\Delta$. Finally, the same holds for the term for $k=0$ and $l=1$,
because of (\ref{eq:majGamma4}), and this proves (\ref{eq:Isa_Y=CPoisson}).

\subsection{Proof of Theorem \ref{theorem:S+S}.}

In the setting of Theorem \ref{theorem:S+S} the measure $G_{\Delta}$ admits
the density $y\mapsto h_{\alpha}(y/\Delta^{1/\alpha})/\Delta^{1/\alpha}$. For
simplicity we set
\begin{equation}
u=\Delta^{\frac{1}{\alpha}-\frac{1}{\beta}}, \label{eq:u}%
\end{equation}
(so $u\rightarrow0$ as $\Delta\rightarrow0$), and a change of variable allows
to write (\ref{eq:pdens_withY}) as
\[
p_{\Delta}(x|\sigma,\beta,\theta,G)={\frac{1}{\theta\Delta^{1/\alpha}}}\int
h_{\beta}\left(  {\frac{x}{\sigma\Delta^{1/\beta}}}-y\right)  h_{\alpha
}\left(  \frac{\sigma y}{\theta u}\right)  .
\]
Therefore
\[
\partial_{\theta}p_{\Delta}(x|\sigma,\beta,\theta,G)=-{\frac{1}{\theta
^{2}\Delta^{1/\alpha}}}\int h_{\beta}\left(  {\frac{x}{\sigma\Delta^{1/\beta}%
}}-y\right)  \breve{h}_{\alpha}\left(  \frac{\sigma y}{\theta u}\right)  ,
\]
and another change of variable in (\ref{eq:I1}) leads to
\begin{equation}
I_{\Delta}^{\theta\theta}=\frac{\theta^{2\alpha-2}u^{2\alpha}}{\sigma
^{2\alpha}}~J_{u}, \label{eq:I1S}%
\end{equation}
where
\begin{equation}
\left\{
\begin{array}
[c]{l}%
J_{u}=\int\frac{R_{u}(x)^{2}}{S_{u}(x)}~dx,\\[2mm]%
R_{u}(x)=\left(  \frac{\sigma}{\theta u}\right)  ^{1+\alpha}\int h_{\beta
}(x-y)\breve{h}_{\alpha}\left(  \frac{\sigma y}{\theta u}\right)  ~dy=\left(
\frac{\sigma}{\theta u}\right)  ^{\alpha}\int h_{\beta}\left(  x-\frac{y\theta
u}{\sigma}\right)  \breve{h}^{\alpha}(y)~dy,\\[2mm]%
S_{u}(x)=\frac{\sigma}{\theta u}\int h_{\beta}(x-y)h_{\alpha}\left(
\frac{\sigma y}{\theta u}\right)  ~dy=\int h_{\beta}\left(  x-\frac{y\theta
u}{\sigma}\right)  h_{\alpha}(y)~dy.
\end{array}
\right.  \label{eq:rs}%
\end{equation}

Below, we denote by $K$ a constant and by $\phi$ a continuous function on
$R_{+}$ with $\phi(0)=0$, both of them changing from line to line and possibly
depending on the parameters $\alpha,\beta,\sigma,\theta$; if they depend on
another parameter $\eta$ we write them as $K_{\eta}$ or $\phi_{\eta}$.
Recalling (\ref{eq:dh^n}), we have
\begin{equation}
h_{\alpha}(x)\sim\frac{c_{\alpha}}{|x|^{1+\alpha}},\qquad\breve{h}_{\alpha
}(x)\sim-\frac{\alpha c_{\alpha}}{|x|^{1+\alpha}},\qquad\int_{\{|y|>|x|\}}%
h_{\alpha}(y)dy\sim\frac{2c_{\alpha}}{\alpha|x|^{\alpha}},\qquad
\mbox{ as $|x|\to\infty$.} \label{eq:dh}%
\end{equation}
Another application of (\ref{eq:dh^n}) when $\beta<2$ and of the explicit form
of $h_{2}$ gives
\begin{equation}
|y|\leq1\quad\Longrightarrow\quad|h_{\beta}^{\prime\prime}(x-y)|\leq
\overline{h}_{\beta}(x):=\left\{
\begin{array}
[c]{ll}%
Kh_{\beta}(x)\qquad & \mbox{if }~\beta<2\\[2mm]%
K(1+x^{2})e^{-x^{2}/2}\qquad & \mbox{if }~\beta=2.
\end{array}
\right.  \label{eq:der}%
\end{equation}

In order to obtain estimates on $R_{u}$ and $S_{u}$, we split the first two
integrals defining these functions into sums of integrals on the two domains
$\{|y|\leq\eta\}$ and $\{|y|>\eta\}$, for some $\eta\in(0,1]$ to be chosen
later. We have $|h_{\beta}(x-y)-h_{\beta}(x)+h_{\beta}^{\prime}(x)y|\leq
\overline{h}_{\beta}(x)|y^{2}$ as soon as $|y|\leq1$, so the fact that both
$f=\breve{h}_{\alpha}$ and $f=h_{\alpha}$ are even functions gives
\[
\left\vert \int_{\{|y|\leq\eta\}}h_{\beta}(x-y)f\left(  \frac{\sigma y}{\theta
u}\right)  dy-h_{\beta}(x)\int_{\{|y|\leq\eta\}}f\left(  \frac{\sigma
y}{\theta u}\right)  dy\right\vert \leq\overline{h}_{\beta}(x)\int
_{\{|y|\leq\eta\}}y^{2}\left\vert f\left(  \frac{\sigma y}{\theta u}\right)
\right\vert ~dy.
\]
On the one hand we have with $f=\breve{h}_{\alpha}$ or $f=h_{\alpha}$, and in
view of (\ref{eq:dh}):
\[
\int_{\{|y|\leq\eta\}}y^{2}\left\vert f\left(  \frac{\sigma y}{\theta
u}\right)  \right\vert ~dy=\frac{\theta^{3}u^{3}}{\sigma^{3}}\int
_{\{|z|\leq\sigma\eta/\theta u\}}z^{2}|f(z)|~dz\leq Ku^{1+\alpha}%
\eta^{2-\alpha}.
\]
On the other hand, the integrability of $h_{\alpha}$ and $\breve{h}_{\alpha}$
and the fact that $\int\breve{h}_{\alpha}(y)dy=0$ yield
\[
\int_{\{|y|\leq1\}}{h}_{\alpha}\left(  \frac{\sigma y}{\theta u}\right)
dy=\frac{\theta u}{\sigma}\int_{\{|z|\leq\sigma/\theta u\}}{h}_{\alpha
}(z)dz=\frac{\theta u}{\sigma}~(1+\phi(u)),
\]%
\[
\int_{\{|y|\leq\eta\}}\!\breve{h}_{\alpha}\left(  \frac{\sigma y}{\theta
u}\right)  dy=\frac{\theta u}{\sigma}\int_{\{|z|\leq\sigma\eta/\theta
u\}}\!\breve{h}_{\alpha}(z)dz=-\frac{\theta u}{\sigma}\int_{\{|z|>\sigma
\eta/\theta u\}}\!\breve{h}_{\alpha}(z)dz=2c_{\alpha}~\frac{(\theta
u)^{1+\alpha}}{\sigma^{1+\alpha}\eta^{\alpha}}~(1+\phi_{\eta}(u)),
\]
Putting all these facts together yields
\begin{equation}
\left\{
\begin{array}
[c]{l}%
\left\vert \int_{\{|y|\leq1\}}h_{\beta}(x-y){h}_{\alpha}\left(  \frac{\sigma
y}{\theta u}\right)  dy-\frac{\theta u}{\sigma}~h_{\beta}(x)\right\vert \leq
u\phi(u)h_{\beta}(x)+Ku^{1+\alpha}\overline{h}_{\beta}(x),\\[2.5mm]%
\left\vert \int_{\{|y|\leq\eta\}}h_{\beta}(x-y)\breve{h}_{\alpha}\left(
\frac{\sigma y}{\theta u}\right)  dy-2c_{\alpha}~\frac{(\theta u)^{1+\alpha}%
}{\sigma^{1+\alpha}\eta^{\alpha}}~h_{\beta}(x)\right\vert \leq u^{1+\alpha
}\left(  h_{\beta}(x)~\frac{\phi_{\eta}(u)}{\eta^{\alpha}}+K\overline
{h}_{\beta}(x)\eta^{2-\alpha}\right)  .
\end{array}
\right.  \label{eq:ru}%
\end{equation}

For the integrals on $\{|y|>\eta\}$ we observe that by (\ref{eq:dh}) we have
$\left\vert h_{\alpha}(\sigma y/\theta u)-c_{\alpha}(\theta u/\sigma
|y|)^{1+\alpha}\right\vert \leq(\theta u/\sigma y)^{1+\alpha}\phi(u)$, and the
same for $\breve{h}_{\alpha}$ except that $c_{\alpha}$ is substituted with
$-\alpha c_{\alpha}$. Then if
\begin{equation}
D_{\eta}(x)=\int_{\{|y|>\eta\}}h_{\beta}(x-y)~\frac{1}{|y|^{1+\alpha}}~dy,
\label{eq:D}%
\end{equation}
we readily get
\begin{equation}
\left\{
\begin{array}
[c]{l}%
\left\vert \int_{\{|y|>1\}}h_{\beta}(x-y){h}_{\alpha}\left(  \frac{\sigma
y}{\theta u}\right)  dy-c_{\alpha}~\frac{(\theta u)^{1+\alpha}}{\sigma
^{1+\alpha}}~D_{1}(x)\right\vert \leq D_{1}(x)u^{1+\alpha}\phi(u),\\[2.5mm]%
\left\vert \int_{\{|y|>\eta\}}h_{\beta}(x-y)\breve{h}_{\alpha}\left(
\frac{\sigma y}{\theta u}\right)  dy+\alpha c_{\alpha}~\frac{(\theta
u)^{1+\alpha}}{\sigma^{1+\alpha}}~D_{\eta}(x)\right\vert \leq D_{\eta
}(x)u^{1+\alpha}\phi(u).
\end{array}
\right.  \label{eq:rub}%
\end{equation}
At this stage, if we put together (\ref{eq:ru}) and (\ref{eq:rub}), we obtain
\begin{equation}
\left\{
\begin{array}
[c]{l}%
\left\vert S_{u}(x)-\left(  h_{\beta}(x)+c_{\alpha}~\frac{(\theta u)^{\alpha}%
}{\sigma^{\alpha}}~D_{1}(x)\right)  \right\vert \leq\left(  h_{\beta
}(x)+u^{\alpha}D_{1}(x)\right)  \phi(u)+Ku^{\alpha}\overline{h}_{\beta
}(x),\\[2.5mm]%
\left\vert R_{u}(x)-c_{\alpha}~\left(  \frac{2h_{\beta}(x)}{\eta^{\alpha}%
}-\alpha D_{\eta}(x)\right)  \right\vert \leq\left(  \frac{h_{\beta}(x)}%
{\eta^{\alpha}}+D_{\eta}(x)\right)  \phi_{\eta}(u)+K\overline{h}_{\beta
}(x)\eta^{2-\alpha}.
\end{array}
\right.  \label{eq:ruq}%
\end{equation}
\smallskip

Our next step is to study the behavior at infinity of the continuous bounded
and positive function $D_{\eta}$. We split the integral in (\ref{eq:D}) into
the sum of the integrals, say $D_{\eta}^{(1)}$ and $D_{\eta}^{(2)}$, over the
two domains $\{|y-x|\leq|x|^{\gamma}\}$ and $\{|y|>\eta,|y-x|>|x|^{\gamma}\}$,
where $\gamma=4/5$ if $\beta=2$ and $\gamma\in((1+\alpha)/(1+\beta),1)$ if
$\beta<2$. On the one hand, $D_{\eta}^{(2)}(x)\leq Kh_{\beta}(|x|^{\gamma})$,
so with our choice of $\gamma$ we obviously have $|x|^{1+\alpha}D_{\eta}%
^{(2)}(x)\rightarrow0$. On the other hand $|x|^{1+\alpha}D_{\eta}^{(1)}(x)$ is
clearly equivalent, as $|x|\rightarrow\infty$, to $\int_{\{|y-x|\leq
|x|^{\gamma}\}}h_{\beta}(x-y)dy$, which equals $\int_{\{|z|\leq|x|^{\gamma}%
\}}h_{\beta}(z)dz$, which in turns goes to $1$. Hence we get for all $\eta
>0$:
\begin{equation}
D_{\eta}(x)~\sim~\frac{1}{|x|^{1+\alpha}}\qquad\mbox{as }~|x|\rightarrow
\infty. \label{eq:equivD}%
\end{equation}

At this stage, we can obtain the behavior of $R_{u}$ and $S_{u}$ as
$u\rightarrow0$. First, an application of (\ref{eq:dh}) to the last formula in
(\ref{eq:rs}) and Lebesgue theorem readily give
\begin{equation}
\lim_{u\rightarrow0}~S_{u}(x)=h_{\beta}(x). \label{eq:CONVS}%
\end{equation}
Also, by (\ref{eq:dh}), (\ref{eq:der}), (\ref{eq:ruq}) and (\ref{eq:equivD}),
we get
\begin{equation}
S_{u}(x)\geq\left\{
\begin{array}
[c]{ll}%
\frac{C}{1+|x|^{1+\beta}} & \mbox{if }~\beta<2\\[2mm]%
C\left(  e^{-x^{2}/2}+\frac{u^{\alpha}}{1+|x|^{1+\alpha}}\right)  \qquad &
\mbox{if }~\beta=2
\end{array}
\right.  \label{eq:CONVSS}%
\end{equation}
for some $C>0$ depending on the parameters, and all $u$ small enough. In the
same way, we see that $u^{\alpha}R_{u}(x)\rightarrow0$, but this not enough.
However, if $r_{\eta}=\frac{2h_{\beta}}{\eta^{\alpha}}-D_{\eta}$ we deduce
from (\ref{eq:ruq}) that for all $\eta\in(0,1]$,
\begin{equation}
\limsup_{u\rightarrow0}\left\vert R_{u}(x)-c_{\alpha}r_{\eta}(x)\right\vert
\leq K\overline{h}_{\beta}(x)\eta^{2-\alpha}. \label{eq:CONVR}%
\end{equation}
A simple computation and the second order Taylor expansion with integral
remainder for $h_{\beta}$ yield
\[
r_{\eta}(x)=\alpha\int_{\{|y|>\eta\}}\frac{h_{\beta}(x)-h_{\beta}%
(x-y)}{|y|^{1+\alpha}}~dy=-\alpha\int_{\{|y|>\eta\}}|y|^{1-\alpha}dy\int
_{0}^{1}(1-v)h_{\beta}^{\prime\prime}(x-yv)dv.
\]
By Lebesgue theorem, $r_{\eta}$ converges pointwise as $\eta\rightarrow0$ to
the function $r$ given by
\[
r(x)=-\alpha\int_{\mathcal{R}}|y|^{1-\alpha}dy\int_{0}^{1}(1-v)h_{\beta
}^{\prime\prime}(x-yv)dv.
\]
Then, taking into account (\ref{eq:CONVR}), and using once more (\ref{eq:ruq})
together with (\ref{eq:dh}), (\ref{eq:der}) and (\ref{eq:equivD}) and also
$\alpha<\beta$, we get
\begin{equation}
\lim_{u\rightarrow0}~R_{u}(x)=c_{\alpha}~r(x),\qquad|R_{u}(x)|\leq\frac
{K}{1+|x|^{1+\alpha}}. \label{eq:CONVRR}%
\end{equation}
\smallskip

We are now in a position to prove (\ref{eq:SS2}), so $\beta<2$ and
$\alpha>\beta/2$ (and of course $\alpha<\beta$). Since $\alpha>\beta/2$, we
see that (\ref{eq:CONVS}), (\ref{eq:CONVSS}) and (\ref{eq:CONVRR}) allow to
apply Lebesgue theorem in the definition (\ref{eq:rs}) to get that
$J_{u}\rightarrow~\int\frac{(c_{\alpha}r(x))^{2}}{h_{\beta}(x)}~dx$: this is
(\ref{eq:SS2}) (obviously $|r(x)|\leq K/(1+|x|^{1+\alpha})$, while $h_{\beta
}(x)>C/(1+|x|^{1+\beta})$ for some $C>0$, so the integral in (\ref{eq:SS2})
converges). \medskip

The other cases are a bit more involved, because Lebesgue theorem does not
apply and we will see that $J_{u}$ goes to infinity. First, we introduce the
following functions:
\[
R^{\prime}(x)=\frac{\alpha c_{\alpha}}{|x|^{1+\alpha}},\qquad S_{u}^{\prime
}(x)=\left\{
\begin{array}
[c]{ll}%
\frac{c_{\beta}}{|x|^{1+\beta}}+\frac{c_{\alpha}\theta^{\alpha}u^{\alpha}%
}{\sigma^{\alpha}|x|^{1+\alpha}}\qquad & \mbox{if }~\beta<2\\[2mm]%
\frac{e^{-x^{2}/2}}{\sqrt{2\pi}}+\frac{c_{\alpha}\theta^{\alpha}u^{\alpha}%
}{\sigma^{\alpha}|x|^{1+\alpha}}\qquad & \mbox{if }~\beta=2
\end{array}
\right.
\]
Below we denote by $\psi(u,\Gamma)$ for $u\in(0,1]$ and $\Gamma\geq1$ the sum
$\phi^{\prime}(u)+\phi^{\prime\prime}(1/\Gamma)$ for any two functions
$\phi^{\prime}$ and $\phi^{\prime\prime}$ like $\phi$ above (changing from
line to line). We deduce from (\ref{eq:dh}), (\ref{eq:der}), (\ref{eq:ruq})
for $\eta=1$, and (\ref{eq:equivD}), that
\begin{equation}
|x|>\Gamma\quad\Rightarrow\quad|R_{u}(x)+R^{\prime}(x)|\leq\psi(u,\Gamma
)R^{\prime}(x)+\left\{
\begin{array}
[c]{ll}%
\frac{K}{|x|^{1+\beta}} & \mbox{if }~\beta<2\\[2mm]%
Kx^{2}e^{-x^{2}/2}\quad & \mbox{if }~\beta=2.
\end{array}
\right.  \label{eq:rrp}%
\end{equation}
In a similar way, we get
\[
|x|>\Gamma\quad\Rightarrow\quad|S_{u}(x)-S_{u}^{\prime}(x)|\leq\left\{
\begin{array}
[c]{ll}%
\psi(u,\Gamma)S_{u}^{\prime}(x) & \mbox{if }~\beta<2\\[2mm]%
\psi(u,\Gamma)S_{u}^{\prime}(x)+K_{0}u^{\alpha}x^{2}e^{-x^{2}/2}\quad &
\mbox{if
}~\beta=2,
\end{array}
\right.
\]
where $K_{0}$ is some constant. Then for any $\varphi>0$ we denote by
$\Gamma_{\varphi}$ the smallest number bigger than $1$, such that $K_{0}%
x^{2}e^{-x^{2}/2}\leq\varphi~\frac{c_{\alpha}\theta^{\alpha}}{\sigma^{\alpha
}|x|^{1+\alpha}}$ for all $|x|>\Gamma_{\varphi}$. The last estimate above for
$\beta=2$ reads as $|S_{u}-S_{u}^{\prime}|\leq S_{u}^{\prime}(\psi
(u,\Gamma)+\varphi)$, so in all cases we have for some fixed function
$\psi_{0}$ as above:
\begin{equation}
S_{u}(x)=S_{u}^{\prime}(x)(1+\rho_{u}(x)),\quad\mbox{where }\quad|\rho
_{u}(x)|\leq\left\{
\begin{array}
[c]{ll}%
\psi_{0}(u,\Gamma) & \mbox{if }~\beta<2,~~|x|>\Gamma\\[2mm]%
\psi_{0}(u,\Gamma)+\varphi\quad & \mbox{if }~\beta=2,~~|x|>\Gamma>\Gamma_{0}.
\end{array}
\right.  \label{eq:rsp}%
\end{equation}

At this stage, we set
\[
J_{u,\Gamma}=\int_{\{|x|>\Gamma\}}\frac{R^{\prime}(x)^{2}}{S_{u}^{\prime}%
(x)}~dx.
\]
Observe that $J_{u}=J_{u,\Gamma}+\sum_{i=1}^{4}J_{u,\Gamma}^{(i)}$, where
\[
J_{u,\Gamma}^{(1)}=\int_{\{|x|\leq\Gamma\}}\frac{R_{u}(x)^{2}}{S_{u}%
(x)}~dx,\qquad J_{u,\Gamma}^{(2)}=\int_{\{|x|>\Gamma\}}\frac{(R_{u}%
(x)+R^{\prime}(x))^{2}}{S_{u}(x)}~dx,
\]%
\[
J_{u,\Gamma}^{(3)}=-2\int_{\{|x|\leq\Gamma\}}\frac{R^{\prime}(x)(R_{u}%
(x)+R^{\prime}(x))}{S_{u}(x)}~dx,\qquad J_{u,\Gamma}^{(4)}=\int_{\{|x|>\Gamma
\}}\frac{R^{\prime}(x)^{2}}{S_{u}(x)}~dx-\int_{\{|x|>\Gamma\}}\frac{R^{\prime
}(x)^{2}}{S_{u}^{\prime}(x)}~dx.
\]

>From (\ref{eq:CONVSS}), (\ref{eq:CONVRR}) and (\ref{eq:rrp}) we get for some
$u_{0}>0$:
\[
\sup_{u\in(0,u_{0}]}J_{u,\Gamma}^{(1)}<\infty,\qquad u\in(0,u_{0}%
]~~\Rightarrow~~J_{u,\Gamma}^{(2)}\leq K+\psi(u,\Gamma)\left(  J_{u,\Gamma
}^{(4)}+J_{u,\Gamma}\right)  .
\]
Cauchy--Schwarz inequality yields
\[
|J_{u,\Gamma}^{(3)}|\leq2\left(  J_{u,\Gamma}^{(2)}\left(  J_{u,\Gamma}%
^{(4)}+J_{u,\Gamma}\right)  \right)  ^{1/2}.
\]
Finally, (\ref{eq:CONVSS}) and (\ref{eq:rsp}) and the definition of
$R^{\prime}$ yield (with $\psi_{0}$ as in (\ref{eq:rsp})):
\[
|J_{u,\Gamma}^{(4)}|\leq\left\{
\begin{array}
[c]{ll}%
2\psi_{0}(u,\Gamma)J_{u,\Gamma} & \mbox{if }~~\beta<2,~~\psi_{0}%
(u,\Gamma)<1/2\\[2mm]%
2(\psi_{0}(u,\Gamma)+\varphi)J_{u,\Gamma} & \mbox{if }~~\beta=2,~~\Gamma
>\Gamma_{\varphi},~~\psi_{0}(u,\Gamma)+\varphi<1/2
\end{array}
\right.
\]

Therefore we get
\[
\left\vert \frac{J_{u}}{J_{u,\Gamma}}-1\right\vert \leq\frac{K}{J_{u,\Gamma}%
}+\psi(u,\Gamma)+2(\psi_{0}(u,\Gamma)+\varphi),
\]
as soon as $\psi_{0}(u,\Gamma)+\varphi<1/2$ and $\Gamma>\Gamma_{\varphi}$, and
with the convention that $\varphi=0$ and $\Gamma_{0}=1$ when $\beta<2$. Then,
remembering that $\lim_{u\rightarrow0}\lim_{\Gamma\rightarrow\infty}%
\psi(u,\Gamma)=0$, and the same for $\psi_{0}$, and that $\varphi=0$ when
$\beta<2$ and $\varphi$ is arbitrarily small when $\beta=2$, we readily deduce
the following fact: Suppose that for some function $u\mapsto\gamma(u)$ going
to $+\infty$ as $u\rightarrow0$, and independent of $\Gamma$, we have proved
that
\begin{equation}
J_{u,\Gamma}\sim\gamma(u)\quad\mbox{as }~u\rightarrow0,\quad\forall\Gamma>1;
\label{eq:EQUI}%
\end{equation}
then $J_{u}\sim\gamma(u)$, and therefore by (\ref{eq:I1S}) we get
\begin{equation}
I_{\Delta}^{\theta\theta}~\sim~\frac{\theta^{2\alpha-2}u^{2\alpha}\gamma
(u)}{\sigma^{2\alpha}}. \label{eq:FINAL}%
\end{equation}
\smallskip

The simplest case is when $\beta<2$ and $\alpha<\beta/2$. Indeed the change of
variables $z=xu^{\alpha/(\beta-\alpha)}$ yields
\begin{align}
J_{u,\Gamma}  &  =\alpha^{2}c_{\alpha}^{2}\int_{\{|x|>\Gamma\}}\frac
{1}{c_{\beta}|x|^{1+2\alpha-\beta}+c_{\alpha}\theta^{\alpha}u^{\alpha
}|x|^{1+\alpha}/\sigma^{\alpha}}~dx\nonumber\\
&  =\frac{\alpha^{2}c_{\alpha}^{2}}{u^{\frac{\alpha(\beta-2\alpha)}%
{\beta-\alpha}}}~\int_{\{|z|>\Gamma u^{\alpha/(\beta-\alpha)}\}}\frac
{1}{c_{\beta}|z|^{1+2\alpha-\beta}+c_{\alpha}\theta^{\alpha}|z|^{1+\alpha
}/\sigma^{\alpha}}~dz,
\end{align}
and we have (\ref{eq:EQUI}) with
\[
\gamma(u)=\frac{\alpha^{2}c_{\alpha}^{2}}{u^{\frac{\alpha(\beta-2\alpha
)}{\beta-\alpha}}}~\int\frac{1}{c_{\beta}|z|^{1+2\alpha-\beta}+c_{\alpha
}\theta^{\alpha}|z|^{1+\alpha}/\sigma^{\alpha}}~dz.
\]
So if we combine this with (\ref{eq:FINAL}), we get (\ref{eq:SS4}). \smallskip

Suppose now that $2\alpha=\beta<2$. Then
\[
J_{u,\Gamma}=2\alpha^{2}c_{\alpha}^{2}\int_{\Gamma}^{\infty}\frac
{1}{x(c_{\beta}+c_{\alpha}\theta^{\alpha}u^{\alpha}x^{\alpha}/\sigma^{\alpha
})}~dx.
\]
For $v>0$ we let $H_{v}$ be the unique point $x>0$ such that $c_{\alpha}%
\theta^{\alpha}v^{\alpha}x^{\alpha}/\sigma^{\alpha}=c_{\beta}$, so in fact
$H_{v}=\rho/v$ for some $\rho>0$. We have
\[
J_{u,\Gamma}\leq\frac{2\alpha^{2}c_{\alpha}^{2}}{c_{\beta}}\int_{\Gamma
}^{H_{u}}\frac{1}{x}~dx+\frac{K}{u^{\alpha}}\int_{H_{u}}^{\infty}\frac
{1}{x^{1+\alpha}}~dx\leq\frac{2\alpha^{2}c_{\alpha}^{2}}{c_{\beta}}%
~\log(1/u)+K.
\]
On the other hand, if $\mu>0$ and $\Gamma<x<H_{u\mu}$ we have $c_{\beta
}+c_{\alpha}\theta^{\alpha}u^{\alpha}x^{\alpha}/\sigma^{\alpha}<c_{\beta
}(1+1/\mu^{\alpha})$, hence
\[
J_{u,\Gamma}>\frac{2\alpha^{2}c_{\alpha}^{2}}{c_{\beta}}\int_{\Gamma
}^{H_{u\varphi}}\frac{1}{x(1+1/\mu^{\alpha})}~dx\geq\frac{2\alpha^{2}%
c_{\alpha}^{2}}{c_{\beta}}~\frac{1}{1+1/\mu^{\alpha}}\log(1/u)-K_{\varphi}.
\]
Putting together these two estimates and choosing $\mu$ big give
(\ref{eq:EQUI}) with $\gamma(u)=\frac{2\alpha^{2}c_{\alpha}^{2}}{c_{\beta}%
}~\log(1/u)$, and we readily deduce (\ref{eq:SS3}). \smallskip

The case $\beta=2$ is treated in pretty much the same way. We have
\[
J_{u,\Gamma}=2\alpha^{2}c_{\alpha}^{2}\int_{\Gamma}^{\infty}\frac
{1}{x^{1+\alpha}(x^{1+\alpha}e^{-x^{2}/2}/\sqrt{2\pi}+c_{\alpha}\theta
^{\alpha}u^{\alpha}/\sigma^{\alpha})}~dx.
\]
We suppose that $\Gamma$ is big enough for $x\mapsto x^{2+2\alpha}%
~e^{-x^{2}/2}$ to be decreasing on $[\Gamma,\infty)$. For $v>0$ small enough,
there is a unique number $H_{v}=x>\Gamma$ such that $c_{\alpha}\theta^{\alpha
}v^{\alpha}/\sigma^{\alpha}=x^{1+\alpha}~e^{-x^{2}/2}/\sqrt{2\pi}$, so in fact
$H_{v}\sim\sqrt{2\alpha\log(1/v)}$ when $v\rightarrow0$. Then
\begin{align}
J_{u,\Gamma}  &  \leq K\int_{\Gamma}^{H_{u}}\frac{e^{x^{2}/2}}{x^{2+2\alpha}%
}~dx+\frac{2\alpha^{2}c_{\alpha}\sigma^{\alpha}}{\theta^{\alpha}u^{\alpha}%
}\int_{H_{u}}^{\infty}\frac{1}{x^{1+\alpha}}~dx\nonumber\\
&  \leq\frac{K}{H_{u}^{\alpha}}+\frac{2\alpha c_{\alpha}\sigma^{\alpha}%
}{\theta^{\alpha}u^{\alpha}H_{u}^{\alpha}}~\leq~\frac{2\alpha c_{\alpha}%
\sigma^{\alpha}}{\theta^{\alpha}u^{\alpha}(2\alpha)^{\alpha/2}~(\log
(1/u))^{\alpha/2}}~(1+o(1)).\nonumber
\end{align}
On the other hand, if $\mu>1$ and $x>H_{u\mu}$ we have $x^{1+\alpha}%
e^{-x^{2}/2}/\sqrt{2\pi}+c_{\alpha}\theta^{\alpha}u^{\alpha}/\sigma^{\alpha
}<+c_{\alpha}\theta^{\alpha}u^{\alpha}(1+\mu^{\alpha})/\sigma^{\alpha}$,
hence
\[
J_{u,\Gamma}>\frac{2\alpha^{2}c_{\alpha}\sigma^{\alpha}}{\theta^{\alpha
}u^{\alpha}}\int_{\Gamma}^{H_{u\varphi}}\frac{1}{x^{1+\alpha}(1+\mu^{\alpha}%
)}~dx\geq\frac{2\alpha c_{\alpha}\sigma^{\alpha}}{\theta^{\alpha}u^{\alpha
}(2\alpha)^{\alpha/2}~(\log(1/u))^{\alpha/2}~(1+\mu^{\alpha})}~(1+o(1)).
\]
So again we see, by choosing $\mu$ close to $1$, that the desired result holds
with
\[
\gamma(u)=\frac{2\alpha c_{\alpha}\sigma^{\alpha}}{\theta^{\alpha}u^{\alpha
}(2\alpha)^{\alpha/2}~(\log(1/u))^{\alpha/2}},
\]
and we deduce (\ref{eq:SS1}).

\bigskip

\bigskip

\bigskip


\bigskip%

\renewcommand{\baselinestretch}{1.1}
\normalsize
\bibliographystyle{ims}

\bigskip%

\Line{\AOSaddress{{Department of Economics}\\
{Princeton University and NBER}\\
{Princeton, NJ 08544-1021}\\
{E-mail: yacine@princeton.edu}}\hfill\AOSaddress
{{Laboratoire de Probabilit\'es (UMR 7599)}\\
{Universit\'e P. et M. Curie (Paris-6)}\\
{75252 Paris C\'edex 05 }\\
{E-mail: jj@ccr.jussieu.fr}}}%

\end{document}